\documentclass{article}
\usepackage{alltt}
\usepackage[all]{xy}
\usepackage{latexsym,amssymb,amsmath,amsfonts, amsthm, xcolor}

\newcommand{\CC}{\mbox{${\rm \:  C\!\!\! I
\;\;}$}}

\def\vbar{\mathchoice{\vrule height6.3ptdepth-.5ptwidth.8pt\kern-.8pt}
  {\vrule height6.3ptdepth-.5ptwidth.8pt\kern-.8pt}
  {\vrule height4.1ptdepth-.35ptwidth.6pt\kern-.6pt}
  {\vrule height3.1ptdepth-.25ptwidth.5pt\kern-.5pt}}
\def\fudge{\mathchoice{}{}{\mkern.5mu}{\mkern.8mu}}
\def\bbc#1#2{{\rm \mkern#2mu\vbar\mkern-#2mu#1}}
\def\bbb#1{{\rm I\mkern-3.5mu #1}}
\def\bba#1#2{{\rm #1\mkern-#2mu\fudge #1}}
\def\bb#1{{\count4=`#1 \advance\count4by-64 \ifcase\count4\or\bba A{11.5}\or
  \bbb B\or\bbc C{5}\or\bbb D\or\bbb E\or\bbb F \or\bbc G{5}\or\bbb H\or
  \bbb I\or\bbc J{3}\or\bbb K\or\bbb L \or\bbb M\or\bbb N\or\bbc O{5} \or
  \bbb P\or\bbc Q{5}\or\bbb R\or\bbc S{4.2}\or\bba T{10.5}\or\bbc U{5}\or
  \bba V{12}\or\bba W{16.5}\or\bba X{11}\or\bba Y{11.7}\or\bba Z{7.5}\fi}}

\newcommand{\NN}{{\bb N}}

\newcommand{\RR}{\mbox{${\rm \:  R\!\!\!\! I
\;\;}$}}

\newcommand{\vs}{\vspace{0.25cm}}

\newtheorem{theorem}{Theorem}
\newtheorem{itlemma}{Lemma}[section]
\newtheorem{itproposition}[itlemma]{Proposition}
\newtheorem{itcorollary}[itlemma]{Corollary}
\newtheorem{itremark}[itlemma]{Remark}
\newtheorem{itremarks}[itlemma]{Remarks}
\newtheorem{itdefinition}[itlemma]{Definition}
\newtheorem{itexample}[itlemma]{Example}

\newenvironment{lemma}{\begin{itlemma}\rm}{\end{itlemma}} 
\newenvironment{remark}{\begin{itremark}\rm}{\end{itremark}} 
\newenvironment{remarks}{\begin{itremarks} \rm}{\end{itremarks}}
\newenvironment{corollary}{\begin{itcorollary}\rm}{\end{itcorollary}}
\newenvironment{proposition}{\begin{itproposition}\rm}{\end{itproposition}}
\newenvironment{definition}{\begin{itdefinition}\rm}{\end{itdefinition}}
\newenvironment{example}{\begin{itexample}\rm}{\end{itexample}}
\newenvironment{fact}{\noindent {{\bf Fact}}:\ \ }{\hfill \medskip}
\newenvironment{claim}{\noindent {\em Claim}. \ \ }{\hfill \medskip}
\newcommand{\be}[1]{\begin{equation}\label{#1}}
\newcommand{\ee}{\end{equation}}
\newcommand{\bl}[1]{\begin{lemma}\label{#1}}
\newcommand{\br}[1]{\begin{remark}\label{#1}}
\newcommand{\brs}[1]{\begin{remarks}\label{#1}}
\newcommand{\bt}[1]{\begin{theorem}\label{#1}}
\newcommand{\bd}[1]{\begin{definition}\label{#1}}
\newcommand{\bp}[1]{\begin{proposition}\label{#1}}
\newcommand{\bc}[1]{\begin{corollary}\label{#1}}
\newcommand{\bfact}[1]{\begin{fact}\label{#1}}
\newcommand{\bex}[1]{\begin{example}\label{#1}}
\newcommand{\ec}{\end{corollary}}
\newcommand{\efact}{\end{fact}}
\newcommand{\eex}{\end{example}}
\newcommand{\el}{\end{lemma}}
\newcommand{\er}{\end{remark}}
\newcommand{\ers}{\end{remarks}}
\newcommand{\et}{\end{theorem}}
\newcommand{\ed}{\end{definition}}
\newcommand{\ep}{\end{proposition}}
\newcommand{\epr}{\end{proof}}
\newcommand{\bpr}{\begin{proof}}
\newcommand{\bcl}{\begin{claim}}
\newcommand{\ecl}{\end{claim}}

\newcommand{\bi}{\begin{itemize}}
\newcommand{\ei}{\end{itemize}}
\newcommand{\ben}{\begin{enumerate}}
\newcommand{\een}{\end{enumerate}}

\usepackage{graphicx}

\title{\bf \Large{On symmetries in time optimal control, sub-Riemannian geometries and the K-P problem}}

\vs

\vs

\author{Francesca Albertini\thanks{Dipartimento di Tecnica e Gestione dei Sistemi Industriali,  Universit\`a di Padova, albertin@math.unipd.it}  \, \, \,  and \, Domenico D'Alessandro\thanks{Department of Mathematics, Iowa State University, Ames, Iowa, U.S.A., e-mail:daless@iastate.edu}}

\begin{document}

\maketitle

\begin{abstract}{
The goal of this paper is to describe a method to solve a class of  time optimal control problems which are   equivalent to finding the sub-Riemannian minimizing geodesics on a manifold $M$. In particular, we assume that the manifold $M$ is acted upon by a group $G$ which is a symmetry group for the dynamics. The action of $G$ on $M$ is  proper but not necessarily free. As a consequence, the {\it orbit space} $M/G$ is not necessarily a manifold but it presents the more general structure of a {\it stratified space}. The main ingredients of the method are a reduction of the problem to the orbit space $M/G$ and an analysis of  the reachable sets on this space. We give general results   relating the stratified structure of the orbit space, and its decomposition into {\it orbit types}, with the optimal synthesis. We consider in more detail the case of the so-called $K-P$ problem where the manifold $M$ is itself a Lie group  and the group $G$ is determined by  a Cartan decomposition 
of $M$. In this case, the geodesics can be explicitly calculated and are analytic. As an illustration, we apply our method and results to the complete optimal synthesis on $SO(3)$
.  }
\end{abstract}

\vs

\vs


\vs

\vs

\section{Introduction} In a recent paper \cite{OptimalityModels}, we have solved the time optimal control problem for a system on $SU(2)$ using a method 
which exploits the symmetries of the problem and provides an explicit description of the reachable sets at every time. In this paper, we formalize such methodology in general and give 
results (proved in sections \ref{SOC} and  \ref{KPPro}) linking the 
structure of a $G$-manifold\footnote{That is, a manifold with the action of a Lie transformation group} to the optimal synthesis. As an example of application,  we provide the complete optimal synthesis  for a minimum time problem on $SO(3)$, which complements some of the results of \cite{BoscaRossi}  obtained with a different method.  

\vs
In order to introduce some of the ideas we shall explore, we provide a brief summary of the treatment of \cite{OptimalityModels} for the problem on $SU(2)$, in its simplest formulation, from the point of view we will take in this paper. The problem is to control in minimum time the system
\be{SU2S}
\dot X=-i u_x \sigma_x X-iu_y \sigma_y X, \qquad X(0)={\bf 1}, 
\ee 
to a desired final condition $X_f \in SU(2)$, subject to a bound on the $L_2$ norm of the control, i.e., $u_x^2+u_y^2 \leq 1$. Here $\sigma_x$ and $\sigma_y$ are {\it Pauli} matrices:
 \be{PauliMat}
\sigma_x:=\begin{pmatrix} 0 & 1 \cr 1 & 0  \end{pmatrix}, \qquad \sigma_y:=\begin{pmatrix} 0 & i \cr -i & 0  \end{pmatrix}, \qquad \sigma_z:=\begin{pmatrix}  1 & 0 \cr 0 & -1\end{pmatrix}. 
\ee
The matrices, $i \sigma_x$ and $i \sigma_y$ span a subspace of 
$su(2)$ which is invariant under the operation of taking a similarity transformation using a {\it diagonal} matrix $D$ in $SU(2)$, i.e., $A \in su(2) \rightarrow DAD^\dagger \in su(2)$. In this respect, the first observation is that if $X:=X(t)$ is an optimal trajectory to go from the identity to $X_f$, then $DXD^\dagger:=DX(t)D^\dagger$ is an optimal trajectory from the identity to $DX_fD^\dagger$. Therefore,  once we have a minimizing geodesic leading to $X_f$, we also  have a minimizing geodesic for every element $DX_fD^\dagger$ in the {\it `orbit'} of $X_f$ and all such geodesics 
project to a unique curve in the space of orbits,  $SU(2)/G$, where $G$ denotes the subgroup of diagonal matrices. The second observation concerns the nature of the orbit space $SU(2)/G$. Since a general matrix in $SU(2)$ can be written as 
\be{genmat}
X:=\begin{pmatrix} x & y \cr -y^* & x^* \end{pmatrix}, \qquad |x|^2+|y|^2=1
\ee
and a similarity transformation by a diagonal matrix 
only affects the phase of the off-diagonal elements, an orbit is uniquely determined by the complex value $x$, with $|x|\leq 1$, i.e., an element of the unit disc in the complex plane which is therefore in one to one correspondence with the elements of $SU(2)/G$. With these facts, we studied in \cite{OptimalityModels} the whole optimal synthesis  in the unit disc. 
Since the problem has a $K-P$ structure (cf. section \ref{KPPro}),  the candidate optimal trajectories can be explicitly expressed in terms of some parameters to be determined according to the desired final condition $X_f$. The number is reduced to only one if we consider the projection on the unit disc of these trajectories. Fixing the time $t$ and 
varying such parameter  we obtained, as parametric  curves, the boundary of the reachable in the unit disc, or, more precisely, the boundary of the projection of the reachable set onto the orbit space. Once an explicit description of the reachable sets is available a method to determine the optimal controls is obtained as a consequence. 

The study of the role of symmetries in optimal control problems is a fundamental  subject in geometric control theory, important both from a conceptual point of view and a practical one as it allows us to reduce the problem to a smaller state (quotient) space. This {\it symmetry reduction} in control problems has a long history (see, e.g., \cite{dodici}, \cite{quattordici}, \cite{quindici}, \cite{sedici}, \cite{venti}, \cite{Martinez}, \cite{trentasei} and see, in particular, \cite{Tomizu} for a recent account). It is obtained  from the application of techniques in {\it geometric mechanics} such as in \cite{Marsden1}, \cite{Marsden2}. However, typically translation of these results of geometric mechanics in control theory has been restricted to the case where the action of the symmetry group $G$ on the underlying manifold $M$ is not only {\it proper} but also {\it free} (definitions are given in section \ref{Backgrou}). In this case the orbit space $M/G$ is guaranteed to be a manifold. In the case where such an action is not free, the orbit space $M/G$ is a {\it stratified space} \cite{LTG}. This is the case discussed here. One example is the above mentioned (closed) unit disc which is a manifold with boundary, a special case of a stratified space.

We have kept the paper as much as possible self contained introducing several concepts from the beginning. In particular,  the paper is organized as follows: In section \ref{Backgrou} we give the necessary background on sub-Riemannian geometry and how it connects with the time-optimal control problem (we refer to \cite{ABB}, \cite{AS} and \cite{Montgomery} for a detailed treatment). This section also contains the basic facts on Lie transformation groups,  in particular the decomposition of the orbit space into {\it orbit types} (see, e.g., \cite{LTG}). In section \ref{SOC}, we present results linking the geometry of the orbit space with the geometry of the optimal synthesis in optimal control. In section  \ref{KPPro}, we apply and expand these results to the case where the problem has an underlying $K-P$ structure. As an example we apply our results to determine the geometry of the optimal synthesis for a control system on $SO(3)$ in section \ref{SO3}. 

\section{Background}\label{Backgrou}

In the next two subsections, we summarize some 
basic concepts in sub-Riemannian geometry and optimal control. We refer to \cite{ABB}, 
\cite{AS},  \cite{Montgomery} for introductory monographs  on  the subject. 

\subsection{Sub-Riemannian structures and minimizing geodesics}\label{SR}

Given a  {\it Riemannian manifold}, $M$, a {\bf sub-Riemannian structure} on $M$ is given by a sub-bundle, $\Delta$, of the tangent bundle $TM$. Letting $\pi_\Delta:  \Delta \rightarrow M$ be the canonical projection, $\Delta$ is a vector bundle on $M$, 
 whose fibers at $x\in M$, $\Delta_x:=\pi_\Delta^{-1}(x) \subseteq T_xM$, are  assumed to have constant  dimension, i.e.,  $\dim{\Delta_x}:=m$ independently  
 of $x$. In the control theoretic setting, a sub-Riemannian structure is often described by giving a set of $m$, linearly independent, smooth vector fields (a {\it frame}) on $M$, 
${\cal F}:=\{X_1,X_2,\ldots,X_m\}$, such that at every point $x \in M$, $\texttt{span}\, \{X_1(x),X_2(x),\ldots,X_m(x)\}=\Delta_x$. It is assumed that ${\cal F}$ is {\it bracket generating}: the smallest Lie algebra of vector fields containing ${\cal F}$, i.e., the Lie algebra {\it generated} by ${\cal F}$, $Lie \, {\cal F}$, is such that, at every point $x \in M$, 
$Lie  \, {\cal F}(x)=T_xM$.  Since $M$ is a  Riemannian manifold, by restricting the Riemannian metric to $\Delta_x \subseteq T_xM$ at every $x \in M$, we obtain a smoothly varying positive definite inner product for vectors in $\Delta_x$, which we will denote by 
$\langle \cdot, \cdot \rangle$. We shall assume that the given frame ${\cal F}$ is {\it orthonormal} with respect to this inner product, that is, $\langle X_j(x), X_k(x) \rangle =\delta_{j,k}$, for every $x \in M$.  

We shall consider {\bf horizontal curves} on $M$. A curve $\gamma: [0,T] \rightarrow M$ is assumed to be {\it Lipschitz continuous} and therefore differentiable almost everywhere in $[0,T]$, with $\dot{\gamma}$ {\it essentially bounded}. That is:  there exists a constant $N$ and  a map $H:[0,T] \rightarrow TM$, with $H(t) \in T_{\gamma(t)}M$, such that 
$\langle H(t), H(t) \rangle_R \leq N$,  for every $t \in [0,T]$, and such that $H(t)=\dot \gamma(t)$, almost everywhere in $[0,T]$. Here $\langle \cdot, \cdot \rangle_R$ denotes the original Riemannian metric on $M$ from which the sub-Riemannian metric $\langle \cdot, \cdot \rangle$ is derived. We shall assume a curve $\gamma$ to be {\it regular}, that is $\dot \gamma(t) \not=0$, almost everywhere in $[0,T]$.  A curve $\gamma$ is said to be {\it horizontal} if $\dot \gamma(t) \in \Delta_{\gamma(t)}$ almost everywhere in $[0,T]$. Given the orthonormal frame ${\cal F}:=\{ X_1,\ldots,X_m\}$, this implies that we can write, almost everywhere in $[0,T]$, 
\be{system}
\dot{\gamma}(t)=\sum_{j=1}^m u_j(t) X_j(\gamma(t)),
\ee
with the functions $u_j$, $j=1,\ldots,m$, given by 
$u_j(t)=\langle X_j(\gamma(t)), \dot{\gamma}(t) \rangle$. We remark that,  because of the smoothness of the $X_j$'s, the continuity of $\gamma$ on the compact set $[0,T]$ and the fact that $\dot \gamma$ is essentially bounded, the functions $u_j$ are also essentially bounded. Therefore a horizontal curve determines $m$ essentially bounded {\it `control'} functions, $u_1,\ldots,u_m$, satisfying (\ref{system}) while, viceversa, given $m$ essentially bounded control functions $u_1,\ldots,u_m$, the solution of (\ref{system}) gives a horizontal curve.   

A horizontal curve $\gamma$  has  a {\bf length}, $l(\gamma)$,  which is given by its 
length in the Riemannian geometry sense, i.e., (using (\ref{system})) 
\be{Riemanniangeolength}
l(\gamma):=\int_0^T \sqrt{\langle \dot{\gamma}(t), \dot{\gamma}(t) \rangle} dt=
\int_0^T\sqrt{\sum_{j=1}^m u^2_j(t)}dt. 
\ee 
A horizontal curve $\gamma$, in the interval $[0,T]$, is said to be {\it parametrized by a constant} if $\langle \dot \gamma, \dot \gamma \rangle $ is constant, almost everywhere in $[0,T]$. It is said {\it parametrized by arclength} if such a constant is equal to one. The image of a curve $\gamma$ in $M$ as well as its length do not change if we {\it re-parametrize} the time $t$. A reparametrization is a Lipschitz, monotone 
and surjective map  $\phi:[0,T^{'}] \rightarrow [0,T]$, and a 
reparametrization of  a curve $\gamma$ is a curve 
$\gamma_\phi :=\gamma \circ \phi \, : \, [0,T^{'}] \rightarrow M$. Given 
a horizontal curve $\gamma$ of length $L$ and $\alpha >0$, consider the increasing map 
$s:[0,T] \rightarrow [0, \alpha L]$, 
\be{reparam}
s(t):=\int_0^t \alpha \|\dot \gamma (r) \| dr,  
\ee 
which is invertible. Let $\phi$ be the inverse map $\phi: [0, \alpha L] \rightarrow [0,T]$. 
Then a standard chain rule argument shows that the re-parametrization $\gamma_\phi:=\gamma \circ \phi$ is parametrized by a constant $\frac{1}{\alpha}$, and in particular it is parametrized by arclength if $\alpha=1$. Viceversa every horizontal curve is the reparametrization of a curve parametrized 
by a constant. We refer to \cite{ABB} (Lemma 3.14 and Lemma 3.15) for details. 

Given two points, $q_0$ and $q_1$, the {\it sub-Riemannian distance} between them,  $d(q_0,q_1)$, is defined as the infimum of the lengths of all horizontal curves $\gamma$, such that 
$\gamma(0)=q_0$, and $\gamma(T)=q_1$. This is obviously greater or equal than the Riemannian distance between the two points where the infimum is taken among all the Lipschitz continuous curves, not necessarily horizontal. The {\bf Chow-Raschevskii theorem} states that if $M$ is connected, in the above described situation and in particular under the bracket generating assumption for ${\cal F}$, $(M,d)$ is a {\it metric space} and its topology as a metric space is equivalent to the one of $M$. This theorem has several consequences including the fact that,  for any two points $q_0$ and $q_1$ in $M$, the distance $d(q_0, q_1)$ is finite, i.e., there exists a horizontal curve $\gamma$ joining $q_0$ and $q_1$ having finite length. Moreover, once $q_0$ is fixed $d(q_0, q_1)$ is continuous as a function of $q_1$. A {\bf minimizing geodesic} $\gamma$ joining $q_0$ and $q_1$, is a horizontal curve  which realizes the sub-Riemannian distance $d(q_0,q_1)$. The {\it existence theorem} says that if $M$ is a complete metric space, and in particular if it is compact, then there exists a minimizing geodesic for any pair of points $q_0$ and $q_1$ in $M$. We shall assume this to be the case in the following.

\subsection{Time optimal control}

The problem we shall consider will be, once $q_0 \in M$ is fixed,  to characterize the minimizing geodesic connecting $q_0$ to $q_1$ for any $q_1 \in M$. This problem is related  to the {\bf minimum time optimal control problem} as described in the following theorem  (cf., e.g., \cite{ABB}).   

\bt{ConnectionOC}
The following two facts are equivalent: 

\begin{enumerate}
\item $\gamma: [0,T] \rightarrow M$ is a minimizing sub-Riemannian geodesic joining $q_0$ and $q_1$, parametrized by constant speed $L$. 

\item $\gamma: [0,T] \rightarrow M$ is a minimum time trajectory of (\ref{system}), subject to $\gamma(0)=q_0$ and $\gamma(T)=q_1$, and subject to $\| \vec u \| \leq L$, almost everywhere. 

\end{enumerate}

\et 

\bpr The proof that $1 \rightarrow 2$ is obtained by contradiction. If $1$ is true and $2$ is not true,  then there exists an essentially bounded conrol $\vec{\tilde{u}} $, with  
$ \|\vec{\tilde{u}} \| \leq L$, and a corresponding solution of (\ref{system}), $\tilde \gamma$, with $\tilde \gamma(0)=q_0$, and $\tilde \gamma(T_1)=q_1$, and $T_1 < T$. Calculate the length of $\tilde \gamma$, 
\be{calculen}
\int_0^{T_1}\| \dot{\tilde{\gamma}} \|dt=\int_0^{T_1} \| \vec{\tilde u} \|dt \leq LT_1 < LT=
\int_0^{T}\| \dot{\gamma} \|dt, 
\ee
which contradicts the fact that $\gamma$ is a minimizing geodesic. 

Let us prove now that $2 \rightarrow 1$. First observe that $\gamma$ must be indeed parametrized
 by constant speed. Since the vector fields $X_j$ in (\ref{system}) are orthonormal,  we know that $\| \dot \gamma \|=\| \vec u\|$ almost everywhere. However 
$\| \vec u\|$ (and therefore $\| \dot \gamma \|$) must be equal to 
$L$ almost everywhere. In fact, assume  $\| \vec u \| < L-\epsilon$, for some $\epsilon >0$,  on an 
interval of positive measure $[t_1, t_2]$, with $\gamma(t_1):=\bar q_1$, and $\gamma(t_2):=\bar q_2$.  Direct computation shows that with the `re-scaled' control 
$\vec{u}_{R}:= \frac{L}{L-\epsilon} \vec u \left(\frac{L}{L-\epsilon} t - 
\frac{\epsilon}{L-\epsilon} t_1 \right)$, the curve $\gamma_R:=\gamma\left(\frac{L}{L-\epsilon}t -\frac{\epsilon t_1}{L-\epsilon} \right)$ is solution of (\ref{system}) with 
$\gamma_R(t_1)=\gamma(t_1)=\bar q_1$ and $\gamma_R(t_1+\frac{L-\epsilon}{L}(t_2-t_1))=\gamma(t_2)=\bar q_2$. Therefore $\vec u_R$, which is an admissible control since its norm  is bounded by $L$, achieves the transfer from $\bar q_1$ to $\bar q_2$ in time $\frac{L-\epsilon}{L} (t_2 - t_1) < (t_2 - t_1)$, which contradicts the optimality of $\vec u$. Moreover $\gamma$ has to be a minimizing geodesic with constant speed $L$. If there was another geodesic $\tilde {\gamma}$ with constant speed $L$, its length would be $LT_1$ which must  be less than the length of $\gamma$, that is $LT$. This implies $T_1 < T$ and contradicts the optimality of the time $T$. 

\epr

In the following we shall assume that our initial point $q_0$ is fixed and we shall look for the sub-Riemannian minimizing geodesics parametrized by constant speed $L$, or equivalently the minimum time trajectories (cf. Theorem \ref{ConnectionOC}) connecting $q_0$ to $q_1$, for any $q_1 \in M$. These curves describe  the so called {\bf optimal synthesis} on $M$. Two loci are important  in the description of the optimal synthesis: The {\bf critical locus} ${\cal CR}(M)$ is the set of points in $M$ where minimizing geodesics loose their optimality, i.e., $p \in {\cal C R}(M)$ if and only if there exists a horizontal curve defined in $[0,T+\epsilon)$, with $T>0 $ and $\epsilon >0$, such that $\gamma(0)=q_0$, $\gamma(T)=p$, $\gamma$ is a minimizing geodesic joining $q_0$ and $\gamma(t)$, for every $t$ in $(0, T)$ and it is not a minimizing geodesic for $t \in (T, T+\epsilon)$. The {\bf cut locus} is the set of points $p \in M$ which are reached by two or more minimizing geodesics, i.e., $p \in {\cal CL}(M)$ if and only if there exists two horizontal curves $\gamma_1$ and 
$\gamma_2$, $[0,T]:\rightarrow M$ such that both $\gamma_1$ and $\gamma_2$ are optimal in $[0,T)$. Because of the existence of a minimizing geodesic, if $p \in {\cal{CL}} (M)$, at least one of the curves $\gamma_1$ and $\gamma_2$ is optimal for $p$, at time $T$. Points in the cut locus are called {\it cut points}. Regularity of minimizing geodesics (cf. \cite{MY}) has consequences on the cut and critical locus. 
Next proposition proves that cut points are also critical points when analyticity is verified, this holds in the $K-P$ problem that will be treated in section \ref{KPPro}. 
We have:

\bp{Propo1}
Assume that all minimizing geodesics are analytic functions of $t$ defined in $[0, \infty)$. Then 
$$
{\cal CL}(M) \subseteq {\cal CR} (M). 
$$
\ep  
\bpr
Assume $p \in {\cal CL}(M)$. Then beside the minimizing geodesic for $p$, $\gamma_1:[0,\infty) \rightarrow M$, with $\gamma_1(T)=p$,  there exists another horizontal curve  $\gamma_2:[0,\infty) \rightarrow M$, which is optimal on $[0,T)$ and satisfies $\gamma_2(T)=p$. At least one between $\gamma_1$ and $\gamma_2$ has to loose optimality at $p$. Therefore,  $p \in {\cal CR} (M)$. If this is not the case,  the concatenation of one of them until time $T$ and the other after time $T$ will  also be optimal, which contradicts analyticity of the minimizing geodesics.  
\epr
 
We shall also consider the {\bf reachable sets} for system (\ref{system}), with $\| \vec u \| \leq L$. The reachable set ${\cal R}(T)$ is the set of all points $p \in M$ such that there exists  an essentially bounded function $\vec{u}$, with $\| \vec u \| \leq L$, a.e., such that the corresponding solution of (\ref{system}), $\gamma$, satisfies, $\gamma(0)=q_0$ and $\gamma(T)=p$. We have that $T_1 \leq T_2$ implies ${\cal R}(T_1) \subseteq {\cal R}(T_2)$. Moreover if $\gamma=\gamma(t)$ is a time optimal trajectory on $[0,T]$ with $\gamma(T)=p$ then $p$ belongs to the boundary of the reachable set ${\cal R}(T)$, which is a closed set 
since the set of values for the control is closed (cf. \cite{Filippov})  
 

\subsection{Symmetries and Lie transformation groups}\label{Symmet}

In addition to the above sub-Riemannian structure, on the manifold $M$, we shall consider the {\bf action of a Lie transformation group} $G$  assuming that it is a {\it left} action,\footnote{i.e., for any $p \in M$, $g_1$ and $g_2$ in $G$, $(g_2 g_1)p=g_2 (g_1 p)$.Every aspect of the theory goes through for right actions with minor modification, that is,  
$g_2(g_1 p)=(g_1 g_2)p$.} it is a {\it proper} action.\footnote{that is, the action map $\alpha: G \times M \rightarrow M \times M$ defined by $\alpha(g,p)=(gp,p)$ is proper, that is, the preimage of any compact set is compact.} We shall also assume that the action map is smooth. We shall denote by $M/G$ the {\it orbit space} of $M$ under the action of $G$, i.e., the space of equivalence classes (orbits) $[p]$, where $p_1$ is equivalent to $p_2$, if and only if there exists a $g \in G$ such that $gp_1=p_2$. $\pi: M \rightarrow M/G$ denotes the canonical projection, and $M/G$ is endowed with the quotient topology. The study of the structure of $M/G$ is part of the {\bf theory of Lie transformation groups}. We now recall the main facts which are needed for our treatment. Details can be found in introductory monographs on the subject, such as, e.g., \cite{LTG}.

Two points $x$ and $y$ in $M$ are said to be {\it of the same type} if their {\it isotropy groups} in $G$ are conjugate. Recall, that the isotropy group of a point $x \in M$, $G_x$,  is the subgroup of elements $g$ of $G$, such that $gx=x$. Two  subgroups $H_1$ and $H_2$ are conjugate if there exists a $g \in G$ such that the map $H_1\rightarrow H_2$, $h \rightarrow g h g^{-1}$ is a group  isomorphism. For any subgroup $H$ of $G$, we denote by $(H)$ the set of groups conjugate to $H$. The subset $M_{(H)} \subseteq M $ is the set of points of $M$ whose isotropy group belongs to $(H)$, or, in other words, the set of points whose isotropy group is conjugate to $H$. There will be only certain classes of groups ${(H)}$ for which $M_{(H)}$ is not empty. These are called the {\bf isotropy types}. It 
is known that $M_{(H)}$ is a {\it submanifold} of $M$ (see, e.g., \cite{Michor1}, 7.4). If two points $x$ and $y$ in $M$ are on the same orbit, i.e., $y=gx$ and $h \in G_y$, then $hy=hgx=y=gx$, so that, $g^{-1}h g \in G_x$. This means that $G_x$ and $G_y$ are conjugate, and therefore $x$ and $y$ both belong to $M_{(G_x)}=M_{(G_y)}$. A consequence of this is that $M_{(H)}$ is the inverse image of a set in $M/G$, $\pi(M_{(H)})=M_{(H)}/G$, which is called the {\it isotropy stratum} of type $(H)$. Isotropy strata have a smooth structure in $M/G$: They are smooth manifolds and the inclusion $M_{(H)}/G\rightarrow M/G$ is smooth (cf.  e.g., \cite{Dimitry}). We remark that $M/G$ itself is not in general a smooth manifold.  It is a smooth manifold if the action of $G$ on $M$ is {\it free} that is the isotropy group $G_x$ is the trivial one given by the identity, for any $x \in M$. In that case, there exists only one possible $(H)$ which contains only the trivial group composed of  only the identity. Therefore   $M_{(H)}/G=M/G$ is a smooth manifold according to the above cited result. In general both $M$ and $M/G$ have the 
structure of a {\bf stratified space.} 
\bd{Stratspace}
A {\bf stratification} of a topological space $N$, is a partition of $N$ into connected manifolds $N_i$, i.e., $N=\cup_{i} N_i$ which is {\it locally finite}, i.e., every compact set in $N$ intersects only a finite number of $N_i$'s. Moreover such a partition satisfies the {\it frontier condition}: If $N_i \cap  \bar N_j \not=\emptyset$ then $N_i \subseteq \bar N_j$ and $\dim(N_i)<\dim(N_j)$.\footnote{The intuitive idea of the frontier connection is that smaller dimensional manifolds in the partition are either totally detached from higher dimensional manifolds (that is the intersection with the closure is empty) or they are part of the boundary.}  
\ed
Consider $M:=\cup_{(H)} M_{(H)}$, where the union is taken over all the isotropy types and further decompose each $M_{(H)}$ into its connected components, so as to obtain a partition of $M$, $M=\cup_{i} M_i$. Moreover, 
partition $M/G$ as $M/G=\cup_i \pi(M_i):=\cup_{i} M_i/G$. Such partitions give a stratification of $M$ and $M/G$, respectively (see, e.g.,\cite{EM} Theorem 1.30). 

On the sets of isotropy types $(H)$ a {\bf partial order} is established by saying that $(H_1) \leq (H_2)$ if $H_1$ is conjugate to a subgroup of $H_2$. This defines subsets 
in $M$ and $M/G$: $M_{\leq H}$ is defined as 
\be{MlessH}
M_{\leq H} :=\cup_{(H_1) \leq (H)} M_{(H_1)}, 
\ee 
with $M_{\leq (H)}/G=\cup_{(H_1) \leq (H)} M_{(H_1)}/G$. We have from the definition that $(H_1)\leq (H_2)$ implies $M_{\leq (H_1)} \subseteq M_{\leq (H_2)}$ and $M_{\leq (H_1)}/G \subseteq M_{\leq (H_2)}/G$.  One of the fundamental results of the theory of transformation groups is the {\bf theorem of existence of minimal orbit type}: {\it There exists a unique orbit type $(K)$ such that $(K) \leq (H)$ for every orbit type 
$(H)$. Moreover $M_{(K)}/G$ is a connected, locally connected, open and dense set in $M/G$, which is a manifold of dimension 
$
\dim{M_{(K)}/G}=\dim{M}-\dim{G}+\dim{K}$ (cf., e.g., \cite{Michor1}) . }
Notice in particular that if $K$ is a discrete group the dimension of $M_{(K)}/G$ is 
$\dim{M}-\dim{G}$. The manifold $M_{(K)}/G$ ($M_{(K)}$) is called the {\bf regular part} of $M/G$ ($M$), while $M/G-M_{(K)}/G$ ($M-M_{(K)}$) is called the {\bf singular part}. 

Given the sub-Riemannian (and Riemannian) structure described in subsection \ref{SR}, with the 
initial point $q_0 \in M$, we shall say that the Lie transformation group $G$ is a {\bf group of symmetries} if the following conditions are verified: (Denote by $\Phi_g$ the smooth map on $M$ which gives the action of $G$, $\Phi_g x:=gx$) 
\begin{enumerate}
\item $q_0$ is a fixed point for the action of $G$ on $M$. That is 
\be{fixedpoint}
\Phi_g(q_0)=q_0, \qquad \forall g \in G, 
\ee 
\item If $\Delta$ denotes the distribution which defines the 
sub-Riemannian structure, the action of $G$ satisfies the 
following invariance property, for every $p \in M$, 
\be{invaprop}
\Phi_{g*} \Delta_p= \Delta_{\Phi_g p}. 
\ee
\item $G$ is a group of {\it isometries}  for the sub-Riemannian metric $\langle \cdot, \cdot \rangle$, that is, for every $p \in M$, and $U,V$ in $\Delta_p$
\be{isomet}
\langle U, V \rangle=\langle \Phi_{g*} U, \Phi_{g*} V \rangle. 
\ee 
 \end{enumerate}  
In the next two sections of this paper, we investigate how the optimal synthesis is related to the orbit type decomposition in a sub-Riemannian structure where  the group $G$ is a group of symmetries for such a structure in the  sense above specified.

\section{Symmetries in the time optimal control problem}\label{SOC}

The following propositions clarify  the role of symmetries and the corresponding orbit space decomposition in the optimal synthesis. 

\bp{Propo1b}
Let $q_1$ and $q_2$ be two points in $M$ on the same orbit, i.e. $q_2=gq_1$ for some  $g \in G$. Let $\gamma_1$ be a minimizing sub-Riemannian geodesic parametrized by constant speed $L$ (and therefore a minimum time trajectory for (\ref{system}) subject to $\|\vec u \| \leq L$ (cf. Theorem \ref{ConnectionOC}), with $\gamma_1(0)=q_0$ and $\gamma_1(T)=q_1$. Then   $\gamma_2:=g\gamma_1$ is a minimizing sub-Riemannian geodesic parametrized by constant speed $L$ (and therefore a minimum time trajectory) as well. 
\ep  
\bpr
First notice that because of property (\ref{invaprop}), for almost every $t$ in $[0,T]$, we have 
\be{horizonthality}
\dot \gamma_2(t)=\Phi_{g*} \dot \gamma_1(t) \in \Delta_{g \gamma_1(t)}=\Delta_{\gamma_2(t)}, 
\ee
so that $\gamma_2$ is horizontal. Moreover because 
of property (\ref{isomet}), a.e., 
\be{isometries2}
\langle \dot \gamma_2(t), \dot \gamma_2(t)\rangle=
\langle \Phi_{g*} \dot \gamma_1,  \Phi_{g*} \dot \gamma_1 \rangle =\langle \dot \gamma_1(t), \dot \gamma_1(t)\rangle=L, 
\ee
so that $\gamma_2$ is also parametrized by constant speed $L$ and $l(\gamma_2)=l(\gamma_1)$.  It is the minimum length since a smaller length would contradict the minimality of 
$\gamma_1$. 
\epr
As a consequence of the previous proposition we have that the space $M/G$ is a metric space with the distance $\bar d$ between the two orbits $\bar q_1$ and $\bar q_2$, defined as  
\be{distbeetorb}
\bar d(\bar q_1, \bar q_2):=\inf_{q_1 \in \bar q_1, q_2 \in \bar q_2} d(q_1, q_2), 
\ee
where $d$ is the sub-Riemannian distance (cf. the Chow-Rashevskii Theorem). A geodesic connecting 
two points $\bar q_1$ and $\bar q_2$, in $M/G$ is a curve that achieves such an infimum. 
\bc{cor1}
The distance $\bar d(\bar q_0, \bar q_1)$ is {\em achieved} by $\pi(\gamma)$ where $\gamma$ is a minimizing sub-Riemannian geodesic connecting $q_0$ to any $q_1$, independently of the representative $q_1 \in \bar q_1$.  
\ec 
Therefore the optimal synthesis in $M$ is the inverse image of the optimal synthesis in $M/G$. Furthermore on $M/G$ we can define {\it critical locus}, {\it cut locus} and {\it reachable sets}, in terms of geodesics in exactly the same way we 
defined them on $M$. These sets are the projections of the corresponding sets in $M$. We have:

\bp{Propo2}

\begin{enumerate}

\item \be{E1} {\cal R}(T)=\pi^{-1}(\pi({\cal R}(T)));
\ee

\item \be{E2} {\cal CR}(M)=\pi^{-1}(\pi({\cal CR}(M))), 
\ee

\item \be{E3}{\cal CL}(M)=\pi^{-1}(\pi({\cal CL}(M))).\ee
\end{enumerate}
\ep 
\bpr
Analogously to the proof of Proposition \ref{Propo1b}, if $q_1$ and $q_2$ are on the same orbit and there is a  horizontal curve $\gamma_1$ connecting $q_0$ to $q_1$, with control $\vec u_1$, $\|\vec u_1 \|=\| \dot \gamma_1 \| \leq L$ in time $T$, then the curve $\gamma_2:=g \gamma_1$, for some $g \in G$, corresponds to control 
$\vec u_2$, with $\|\vec u_2 \|=\|\dot \gamma_2 \|=\| \dot \gamma_1\| \leq L$, a.e., connecting $q_0$ to $q_2$, in the same time $T$. Therefore $q_1$ is in ${\cal R}(T)$ if and only if  $q_2$ is in ${\cal R}(T)$, which proves (\ref{E1}). Analogously we can prove that if $q_1$ and $q_2$ are on the same orbit they are both in ${\cal{CL}}(M)$ and ${\cal{CR}}(M)$  or none  of them is, that is, (\ref{E2}) and (\ref{E3}) also hold. We illustrate the proof for ${\cal CL}(M)$ (the proof of ${\cal CR}(M)$ is similar). Assume by contradiction that $q_1 \in {\cal CL}(M)$ and $q_2 \notin {\cal CL}(M)$, with $q_2= g q_1$, for $g \in G$. Since $q_1 \in {\cal CL}(M)$, there are two different minimizing sub-Riemannian geodesics with  $q_1$ as their final point, $\gamma$ and $\tilde \gamma$.  Then $g\gamma$ and  $g\tilde \gamma$ will be two different  geodesics
 with $q_2$ as their final point. In fact, $g \gamma(t)=g \tilde \gamma(t)$, for all $t \in [0,T]$ would imply $\gamma(t)=\tilde\gamma(t)$, which we have excluded.  
\epr 

\vs

\bp{Propo3} Assume that $q_1$ does not belong to the cut locus, i.e., $q_1 \notin {\cal CL}(M)$ and $\gamma$ is a sub-Riemannian minimizing geodesic in $[0,T]$ for $q_1$. Then for every $t \in (0,T)$, we have the following relation for the isotropy groups 
\be{relaisot}
G_{q_1} \subseteq G_{\gamma(t)}. 
\ee
\bpr If (\ref{relaisot}) is not valid then there exists a $g \in G_{q_1}$ with 
$g \notin G_{\gamma (t)}$. Then the curve $g\gamma$ is also a geodesic since it has the same length as $\gamma$ (cf. Proposition \ref{Propo1b}). Moreover since $g \in G_{q_1}$ the geodesic  $g\gamma$ goes to $q_1$,  but it is different from $\gamma$ since $g \notin G_{\gamma (t)}$, for some $t \in (0,T)$. Therefore $q_1$ must belong to the cut locus. Something we have excluded. 
\epr
\ep

Under the assumption that geodesics are analytic we can obtain  for general $q_1 \in M$ the converse inclusion 
to (\ref{relaisot}). 
\bp{Propo7} If all geodesics are analytic then, {\em for any $q_1 \in M$}, and any 
geodesic $\gamma:[0,T] \rightarrow M$,  connecting $q_0$ to $q_1$, we have for any $t\in (0,T)$ 
\be{inverserela}
G_{\gamma(t)} \subseteq G_{q_1}. 
\ee
\ep
\bpr
The proof uses some ideas of Lemma 3.5 in \cite{Dimitry}. Assume there exists a $\bar{t}\in (0,T)$ and a $g\in G_{\gamma(\bar{t})}$ which is not in $G_{q_1}$. Then the curve which is equal to $\gamma$ between $q_0$ and $\gamma(\bar{t})$ and is equal to $g\gamma$ between $\gamma(\bar{t})$ and $g q_1$ 
(which is different from $q_1$), has the same length as 
the curve $g \gamma$. Such a curve is a minimizing geodesic, according to Proposition \ref{Propo1}, since it has the same length as $\gamma$, that is the minimal length. However such a curve which is equal to $\gamma$ in an interval of positive measure until time $t$ and different from $\gamma$ afterwards cannot be analytic, which contradicts the analyticity of all the geodesics.  
\epr

We collect in the following Corollary some consequences of Propositions \ref{Propo3} and  \ref{Propo7}. The Corollary describes how minimizing sub-Riemannian geodesics sit in the orbit type decomposition of $M$ and $M/G$. 

\bc{coro77}
If all minimizing sub-Riemannian 
geodesics are analytic, then for every $q \in M$ ($\bar q \in M/G$), every minimizing geodesic is entirely contained in $M_{\leq ({G_q})}$ ($M_{\leq ({G_q})}/G$) except possibly for the initial point $q_0$. In particular if $q \in M_{(K)}$ (the regular part of $M$)  then the whole minimizing geodesic except possibly the initial point $q_0$ is in  $M_{(K)}$. If, in addition, $q \notin {\cal CL}(M)$ then  
\be{egalite}
G_{\gamma(t)}=G_{q}, \qquad \forall t \in (0,T)
\ee
and the corresponding  sub-Riemannian geodesic is all contained in $M_{(G_q)}$ 
\ec  
\noindent Similar `convexity' results in the Riemannian case are given in (\cite{Dimitry} 3.4).
Corollary \ref{coro77} gives a general principle for the behavior of geodesics in the presence of a group of symmetries for the optimal control problem:

{\it  The geodesics can only go from lower ranked strata such as the lowest $M_{(K)}$ to higher ranked ones but not viceversa. If a geodesic touches a higher ranked point and then goes back to a lower ranked one, it means that it has lost optimality and therefore the point belongs to the critical locus.}

\br{sugge3}
The corollary suggests that the points in the singular part of $M$, $M_{sing}:=M-M_{(K)}$,   are, in general, good candidates to be in the cut locus ${\cal CL}(M)$.\footnote{and therefore in the critical locus ${\cal{CR}}(M)$ cf. Proposition \ref{Propo1}} In fact, any point $q \in M_{sing}$ which has a geodesic with points in $M_{(K)}$, which is an {\it open and dense set} in $M$, 
must be in ${\cal CL}(M)$. Points in  $M_{sing}$ which are {\it not} in ${\cal CL}(M)$ must be such that {\it every geodesic} leading to that point must be entirely contained in the singular part of $M$ since we have that (\ref{egalite}) is verified. In the $SU(2)$ example of \cite{OptimalityModels} it holds that ${\cal CL}(M)={\cal CR}(M)=M_{sing}$. However this 
is not always the case and in general the situation changes by changing  the group of symmetries  we consider (cf. Remark \ref{connesso} below). 
\er

\section{The K-P problem}\label{KPPro}

An example of a subRiemannian problem with symmetries is the {\bf $K-P$ problem} discussed in \cite{OptRes}, \cite{BoscaRossi}.  In this section, we shall focus on this type of problems. 

In the $K-P$ problem, the manifold  $M$ is a {\it semisimple  Lie Group} with its Lie Algebra of right invariant vector fields $\mathcal{L}$. The Lie algebra   $\mathcal{L}$ has   a {\bf Cartan decomposition}, that is, a vector space decomposition 
\be{CDD1}
\mathcal{L}= {\mathcal K} \oplus {\mathcal P},
\ee
with the  commutation relations:\footnote{It is true (see \cite{OptRes}, Appendix A; see also Lemma 3.4 and Corollary 3.5 in \cite{NoiSIAM}) that for  semi-simple Lie Algebras the equality must hold in the second and the third of these inclusions.}
\be{Cartandecodef}
 [{\cal K}, {\cal K}] \subseteq {\cal K}, \quad [{\cal K}, {\cal P}] \subseteq {\cal P}, \qquad [{\cal P}, {\cal P}] \subseteq {\cal K}.
\ee
 The Lie algebra ${\cal L}$ is endowed with a (pseudo)-inner product defined 
 by the {\it Killing form} $\langle A,D \rangle:=Kill(A,D):=Tr(ad_Aad_D)$, where $ad_A$ is the linear operator on ${\cal L}$, given $ad_A(X)=[A,X]$.
{A comprehensive introduction to notions of Lie theory can be found for example in \cite{Knapp}}.  In particular since ${\cal L}$ is assumed semisimple,  
 $Kill$ is non degenerate (this is the Cartan criterion for semisimplicity). Associated with a Cartan decomposition is a {\it Cartan involution}, that is,  an automorphism $\theta$ of ${\cal L}$, such that $\theta^2$ is the identity, and ${\cal K}$ and ${\cal P}$ above are the $+1$ and $-1$ eigenspaces of $\theta$ in ${\cal L}$. Moreover $B( A, D):=-Kill(A,\theta D)$ is a {\it positive definite} bilinear form and therefore an inner product defined on all of ${\cal L}$. Notice that this, in particular, implies that ${\cal K}$ is a compact subalgebra of ${\cal L}$, (i.e. $Kill(A,A)<0$, if $A\in {\cal K}$ and $A \not=0$) and ${\cal K}$ and ${\cal P}$ are orthogonal with respect to such inner product.\footnote{If $A\in {\cal K}$ and $D \in {\cal P}$
 \[
 B(A, D ) :=-Kill(A,\theta D)=Kill(A,D)=Kill(\theta A, \theta D)=
 Kill(A,\theta D)=-B( A, D ), 
 \]  
 where we have used the property of the Killing form that for every Lie algebra 
 automorphism $\phi$, $Kill(A,D)=Kill(\phi A, \phi D)$.} Using the inner product $B$ on the Lie algebra ${\cal L}$,  one naturally defines a Riemannian metric on $M$. In fact, for any point $b \in M$ and any tangent vector $U \in T_bM$, one can associate a right invariant vector field $X_U$ defined as $X_U|_a:=R_{b^{-1}a *}U$ and this association is an isomorphism of vector spaces. Then the Riemannian metric 
 $\langle \cdot, \cdot \rangle_R$ is defined as (with $U,V \in T_bM$)
\be{Riemannianmet}
\langle U, V \rangle_R:=B(X_U, X_V). 
\ee
 The {\bf $K-P$ problem} is the minimum time problem for system (\ref{system}) on a Lie group $M$ with a Cartan decomposition as described above,  
 where the vector fields $X_j$ are right-invariant vector 
 fields\footnote{Notice that we could have as well set 
 up the whole treatment for {\it right} 
 invariant vector fields but we could have given an analogous 
 treatment for {\it left}  invariant vector fields}  on $M$ spanning  ${\cal P}$ and orthonormal with respect the inner product $B$. The initial point $q_0$ is the identity of the group $M$. The problem is to steer from $q_0$ to an arbitrary final condition in $M$ in minimum time,  subject to the condition $\| \vec u \| \leq L$. 
 
 The problem can be cast in the above sub-Riemannian setting with symmetries as follows: The distribution of  vector fields ${\cal P}$ in ${\cal L}$ defines a sub-Riemmannian structure  on $M$ with the sub-Riemannian metric 
 at any point $b$ defined by the restriction of  
 $\langle \cdot, \cdot \rangle_R$ to ${\cal P}|_b$ for any $b$ in $M$. Consider now a Lie subgroup of $M$, $G$ (not necessarily connected) with Lie algebra ${\cal K}$, which acts on $M$ by conjugation, i.e., for $p \in M$, $g \in G$, $\Phi_g(p):=gp:
 = g \times p \times g^{-1}$, where $\times$ is the group operation in $M$. Such (left) action induces a map on the Lie algebra ${\cal L}$ given by its differential $\Phi_{g *}$ which is a Lie algebra automorphism. We assume that the map $\Phi_g$ is an {\it isometry} and that for every connected component $j$ of $G$ there exists a $g_j$ such that $\Phi_{g_j *} {\cal P} \subseteq {\cal P}$. This also implies, because of the (Killing) orthogonality of ${\cal K}$ and ${\cal P}$,  $\Phi_{g_j *} {\cal K} \subseteq {\cal K}$. Moreover, because of (\ref{Cartandecodef}) these properties are not restricted to $g_j$ but are true for every $g \in G$.\footnote{Consider a connected component of $G$. We know that there exists a number of right invariant vector fields $X_1,X_2,\ldots, X_m$ in ${\cal K}$ such that denoting by $\sigma_{1,t}$, $\sigma_{2,t}$,..., $\sigma_{m,t}$ the corresponding flows, we have $\sigma_{m,t_m} \circ \sigma_{{m-1},t_{m-1}} \circ \cdots \circ \sigma_{1,t_1} (g_j)=g$. For every $r=1,\ldots m$ the map $\sigma_{r,t}$ is real analytic as a function of $t$. Denote by $\bar g:=\sigma_{1,t_1}(g_j)$. We want to show that $\Phi_{\bar g * }{\cal P} \subseteq {\cal P}$ and applying this $m$ times we have that $\Phi_{g*} {\cal P} \subseteq {\cal P}$. Consider $K$ in ${\cal K}$ and $P \in {\cal P}$ and the Killing inner product $B(K, \Phi_{\sigma_{1, t}(g_j) *} P)$ which is a real analytic function of $t$ at every point in $M$ and it is zero for $t=0$. By taking the $k$-th derivative of this function at $t=0$, we obtain, using the definitions of Lie derivative 
$$
\frac{d^k}{dt^k}|_{t=0} B(K, \Phi_{\sigma_{1, t}(g_j) *} P)=
B(K, \frac{d^k}{dt^k}|_{t=0}\Phi_{\sigma_{1, t}(g_j) *} P)=
B(K, ad_{X_1}^k\Phi_{g_j *} P)=0, 
$$ 
where $ad_{X_1}^k$ denotes the $k-$the repeated Lie bracket with $X_1$ and we have used (\ref{Cartandecodef}).} A special case is when $G$ is the connected component containing the identity with Lie algebra ${\cal K}$, in which case $g_j$ can be taken equal to the identity.  
 
In the following we shall restrict ourselves to {\bf linear Lie groups} so that $M$ and $G$ will be  Lie groups of matrices.\footnote{The example of $SU(2)$ treated in \cite{OptimalityModels} and the example of $SO(3)$ of the next section are $K-P$ problems of this type.} In the standard coordinates (inherited from the standard ones  of $Gl(n,\RR)$ or $Gl(n,\CC)$) the system (\ref{system}) is written as 
\be{systemMatrices}
\dot X(t) =  \sum_{j} B_j X u_j(t),
\ee
where the Lie algebra ${\cal L}$ of matrices of the Lie group $M$ has a Cartan decomposition as in (\ref{CDD1}) and (\ref{Cartandecodef}) and the $B_j$'s span an orthonormal basis of ${\cal P}$.\footnote{Here with minor abuse of notation we identify ${\cal K}$, ${\cal L}$ and ${\cal P}$ with the spaces of matrices representing the corresponding vector fields.} The $K-P$ problem is the 
minimum time problem, with initial condition $X_0= {\bf 1}$ equal to the identity matrix subject to $||\vec u|| \leq L$. The symmetries are  given  by the transformations $X\rightarrow KXK^{-1}$, for $K \in G$, which induce transformations  on the matrices $B \in {\cal L}$, $B \rightarrow K B K^{-1}$. These are symmetries because they preserve ${\cal P}$ and ${\cal K}$, and the commutation relations (\ref{Cartandecodef}). 

As a special case of what we have seen in general, the minimum time 
control for system (\ref{systemMatrices}) is equivalent to that of finding 
minimizing geodesics on $M$ and it can be treated  on $M/G$. On the orbit space 
we can describe the cut locus, the critical locus and the reachable sets. 
\br{condrift} The knowledge of the reachable sets  for a $K-P$ problem  of the form (\ref{systemMatrices}) also gives the reachable sets for the larger class of systems 
\be{systemMatriceswith}
\dot X =  AX+\sum_{j}^m B_j X u_j(t),
\ee
with the drift $AX$, with $A \in {\cal K}$ and the $B_j \in {\cal P}$. In fact,  
consider the change of coordinates $U(t):=e^{-At} X(t)$. A straightforward 
calculation gives 
$$
\dot U=\sum_{j=1}^m e^{-At} B_j e^{At} U u_j, 
$$
and since $B \rightarrow e^{-At} B e^{At}$ is assumed to be an isometry, there exists an orthogonal matrix $a_{j,k}:=a_{j,k}(t)$ such that  
$ e^{-At} B_j e^{At}=\sum_{k=1}^m a_{j,k}(t) B_k$, so that we have 
\be{system2}
\dot U=\sum_{j=1}^m \sum_{k=1}^m a_{j,k}(t) B_k U u_j=\sum_{k=1}^m B_kUv_k(t),   
\ee
where $v_k(t):=\sum_{j=1}^m a_{j,k}(t)u_j(t)$, and $\|\vec v \| \leq L$ if 
and only if $\| \vec u \| \leq L$. Therefore the reachable set of system   (\ref{system2}) coincides with the one of (\ref{system}) and knowledge of the reachable 
set for system (\ref{system2}), ${\cal R}_U(t)$, gives the reachable set for system 
(\ref{systemMatriceswith}), ${\cal R}_X(t)$, via the relation ${\cal R}_X(t)=e^{At}{\cal R}_U(t)$. 
\er    
 
\vs

In the $K-P$ problem, the equations of the Pontryagin Maximum Principle are {\it explicitly integrable,} and give (cf. \cite{OptRes} and references therein) that the optimal control $\vec u$ is such that there exist matrices $A_k \in {\cal K}$ and $A_p \in {\cal P}$ with 
\be{espressioOptCont}
\sum_{j=1}^m B_j u_j(t)=e^{A_k t} A_p {e^{-A_k t}}, 
\ee
with $\|A_p\|=L$.  Therefore, the  optimal trajectories  satisfy 
\[
\dot X= e^{A_k t} A_p {e^{-A_k t}} X, \quad X(0)={\bf 1},
\quad A_k \in {\cal K}, \quad A_p \in {\cal P}, 
\]
and the solution can be written explicitly as 
\be{geodesicKP}
X(t)=e^{-A_k t}e^{(A_k +A_p)t}.
\ee
The  geodesics (\ref{geodesicKP}) are {\bf analytic} curves. Therefore all 
the results on the geometry of the optimal synthesis in the 
previous section apply. Moreover, for every geodesic in the 
orbit space $M/G$  (which is the  projection of a geodesic in $M$), we can always 
take a representative (\ref{geodesicKP}) in $M$ with $A_p:=A_a \in {\cal A}$, with ${\cal A}$ a {\it maximal Abelian (Cartan) subalgebra} in ${\cal P}$. To see this, recall the  known property of the Cartan decomposition that if 
${\cal A} \subseteq {\cal P}$ is a maximal Abelian subalgebra in ${\cal P}$, then 
\be{Abeliansub}
{\cal P}=\bigcup_{K \in {G}} K {\cal A} K^{-1}.
\ee
Therefore we can write (\ref{espressioOptCont}), for $A_a \in {\cal A}$, as 
\be{AbeliansubA}
\sum_{j=1}^m B_j u_j(t)=Ke^{\bar A_k t} A_a {e^{-\bar A_k t}} K^{-1}, 
\ee
for $K \in G$ and $\bar A_k \in {\cal K}$. 
Using $A_p:=KA_aK$, with $K \in G$, we have (cf. (\ref{geodesicKP})) 
\be{equivclass}
\left[ e^{-A_k t} e^{(A_k+A_p) t} \right]=
\left[ e^{-A_k t} e^{(A_k+KA_a K^{-1}) t} \right] =
\left[ K e^{-\bar A_k t} e^{(\bar A_k+A_a) t} K^{-1} \right ]= \left[ e^{-\bar A_k t} e^{(\bar A_k+A_a) t} \right], 
\ee
 with $\bar A_k:=K^{-1} A_k K$.

The following proposition gives some restrictions 
on the pairs $(A_k,A_p)$ for points that are not on the cut locus of $M$. This proposition can be used to prove that a certain point is in the cut locus.\footnote{This is done for example in the next section in Proposition \ref{seconda}.} 
 
\bp{ProposizioneFrancesca}
Let $X_f \notin {\cal CL}(M)$. Let $H$ denote the isotropy group of $X_f$. Then the pair $(A_k,A_p)$ giving the minimizing geodesic are such that for every $\hat H \in H$, 
\be{invariance}
\hat H ad_{A_k}^n A_p \hat H^{-1}=ad_{A_k}^n A_p,  
\ee
for any $n \geq 0$. 
\ep
\bpr
We know from Propositions \ref{Propo3} and  \ref{Propo7} and from formula (\ref{geodesicKP}) that the pairs $(A_k,A_p)$ satisfy the invariance property  
\be{inv2a}
\hat H e^{-A_k t}e^{(A_k +A_p)t} \hat H^{-1}=e^{-A_k t}e^{(A_k +A_p)t},  
\ee
for every $t \in [0,T]$, where $T$ is the minimum time associated to $X_f$. Taking the $n-$th derivative and by induction it is seen that 
this implies that $X^n(t)$ also satisfies the invariance property with respect to $\hat H$, i.e., 
\be{inva3}
\hat H X^{n}(t) \hat H^{-1}=X^n(t), 
\ee
where $X^n(t)$ is defined as 
\[
X^{(n)}(t):= e^{-A_kt}H_ne^{(A_k+A_p)t}.
\]
with
\be{inductive}
\begin{array}{l}
H_0={\bf 1}, \\
H_{n+1}=H_nA_p+[H_n,A_k].
\end{array}
\ee
Using the invariance (\ref{inva3}) of $X^n(t)$ at $t=0$, it follows that all the matrices $H_n$ are also $\hat H$-invariant. We want to show that this implies the invariance of   
$$
L_n:=ad_{A_k}^{n-1} A_p, 
$$ 
for each $n \geq 1$. For $n=1$ $L_1=A_p$ and it is clear that $L_n$ is invariant since $L_1=H_1=A_p$. From this, we proceed by induction on $n$, for $n \geq 2$. We shall prove that each $H_n$, $n\geq 2$, can be written, { with $l_n=2^{n-1}-1$}, as:
\be{Hn}
H_n=\sum_{j=1}^{l_n} V_j^nW_j^n    +  L_n,
\ee
with $V_j^n$ and $W_j^n$ invariant and both addends of some $H_s$, with $s<n$. From this, since $H_n$ is also invariant, we must have that $L_n$ is invariant as well.

First notice that $H_2=A^2_p+ [A_p,A_k] =A_pA_p +L_2$, so clearly the statement holds for $n=2$.

Assume that the statement holds for $H_n$, then we have:
\be{formulabas}
H_{n+1}= H_nA_p+[H_n,A_k]= \left(\sum_{j=1}^{l_n} V_j^nW_j^nA_p    +  L_nA_p\right) +\left( \sum_{j=1}^{l_n} [V_j^nW_j^n,A_k]    +  [L_n,A_k]\right).
\ee
Letting:
\[
V_j^{n+1}=V_j^nW_j^n, \   \  W_j^{n+1}=A_p,  \   \   V_{l_{n}+1}^{n+1}=L_n,  \   \   W_{l_{n}+1}^{n+1}=A_p.
\]
we can write 
\[
\left(\sum_{j=1}^{l_n} V_j^nW_j^nA_p    +  L_nA_p\right)= \sum_{j=1}^{l_n+1} V_j^{n+1}W_j^{n+1},
\]
where $V_j^{n+1}$ are invariant, since product of invariant, and they are addend of $H_n$, and $W_j^{n+1}=A_p$ is also invariant and it is in $H_1$.
Now we have in (\ref{formulabas})
\[
[V_j^nW_j^n,A_k]= V_j^n[W_j^n,A_k] +[V_j^n,A_k]W_j^n,
\]
If $W_j^n$ is one of the addends of $H_s$ with $s<n$, then $[W_j^n,A_k]$ is one of the addends  of $H_{s+1}$, so it is also invariant, by inductive assumption, and $s+1<n+1$, so 
$V_j^n[W_j^n,A_k]$ is the product of two invariant factors which are addends of two $H_p$ for some $p<n+1$. The same argument applies to
 $[V_j^n,A_k]W_j^n$. So also the sum in the second brackets can be rewritten in the desired form.
 Now it is sufficient to notice that $[L_n,A_k]=L_{n+1}$.
\epr

\subsection{A method to obtain the optimal synthesis for $K-P$ problems}

The previous considerations suggests a general methodology 
to find the optimal synthesis for time optimal control problems with symmetries and in particular for $K-P$ problems. 

The {\bf first} step of the method is to identify a group of symmetries. 
There are in general several choices of groups, connected and not connected. In the $K-P$ case the natural choice is the connected Lie group corresponding to the subalgebra ${\cal K}$ in the Cartan decomposition, or a possible not connected Lie group having ${\cal K}$ as its Lie algebra. It is typically convenient to take the Lie group $G$ as large as 
possible so as to have a finer orbit type decomposition of $M/G$, which we would like to have of as small dimension as possible. 

The {\bf second} step of the procedure is to determine the nature of $M/G$ so that the problem is effectively reduced to a lower dimensional space. This is important both from a conceptual and practical point of view since a computer solution of the problem will have to consider a smaller number of parameters. This task typically requires some analysis since not all the quotient spaces are known in the literature.\footnote{Typical cases in the literature  
look at a Lie group $M$ where the conjugation action on $M$ is 
given by $M$ itself and not by a subgroup $G$ of $M$ as in our case.} An analysis of the various isotropy groups of the points in $M$ reveals the stratified structure of $M/G$ which, as we have seen in section \ref{SOC},  has consequences for the optimal synthesis.

The {\bf third step} is to obtain the boundaries of the reachable 
sets in $M/G$, that is, the projections of the boundaries of 
the reachable sets ${\cal R}(t)$ in  $M$. In order to do this,
if $A_a$ is an element in the Cartan subalgebra ${\cal A} \subseteq {\cal P}$, we write a representative of a geodesic as  (cf. (\ref{equivclass})), 
\be{representative}
X(t):=e^{-\bar A_k t} e^{(\bar A_k + A_a)t}, 
\ee
with $\bar A_k \in {\cal K}$ and $A_a \in {\cal A}$ 
and $\|A_a\|=L$. By fixing $t$ and varying $\bar A_k \in {\cal K}$ and $A_a \in {\cal A}$ 
we obtain an hyper-surface in $M/G$, part of which is the boundary of 
the reachable set at time $t$. The determination of the sets in ${\cal K}$ and 
${\cal A}$ which is mapped to this boundary is an 
analysis problem to be considered on a case by case basis, which is obviously simpler in low dimensional cases, and requires help from computer simulations in higher dimensional cases. 

The {\bf fourth step} is to find the first $t$ such that $\pi({\cal R} (t))$ contains $\pi(X_f)$. At this value of $t$, there are matrices $A_k$ and $A_a$ such that  
$[e^{-\bar A_k t} e^{(\bar A_k + A_a)t}]=\pi(X_f)$. 

Finally the {\bf fifth step} is to find $K \in G$ such that 
\be{finalstep3}
K e^{-\bar A_k t} e^{(\bar A_k + A_a)t} K^{-1}=X_f.  
\ee
This gives the correct pair $(A_k, A_p)$ to be used in the optimal control (\ref{espressioOptCont}): $A_k=K\bar A_k K^{-1}$, $A_p:= K \bar A_a K^{-1}$. 
From the last two steps, it follows that the problem is therefore effectively divided in two. Restricting ourselves to the orbit space we first find an optimal control to drive the state of the system to the desired orbit. Then, in the fifth step, we move inside  the orbit to find exactly the final condition we desire.

The treatment of the optimal synthesis on $SO(3)$ in the following section gives an example of application of this method.    
 
\section{Optimal synthesis for the $K-P$ problem on $SO(3)$}
\label{SO3}

A basis of the Lie algebra  of skew-symmetric real  
$3\times 3$ matrices, $so(3)$, is given by:
\[
p_1:=\begin{pmatrix} 0 & 0 & 0 \cr 
0 & 0 & - 1 \cr
0 & 1 &  0
\end{pmatrix},  \   \   \  p_2:=\begin{pmatrix} 0 & 0 & 1 \cr 
0 & 0 &0  \cr
-1 & 0 &  0
\end{pmatrix},  \  \   \  k:=\begin{pmatrix} 0 & -1 & 0 \cr 
1 & 0 & 0 \cr
0 & 0 &  0
\end{pmatrix}.  
\]
We consider the $K-P$ Cartan decomposition of $so(3)$ where ${\mathcal{K}}=\text{span} \{k\}$, and ${\mathcal{P}}=\text{span} \{ p_1, \, p_2\}$.
 There are two possible maximal groups of symmetries with Lie algebra ${\cal K}$. A {\it maximal connected} Lie group,  ${\bf K}^+$,  which is the connected component containing the identity  and consists of  matrices of the form 
\be{forma1}
K^+(r):= \begin{pmatrix} \cos(r) & \sin(r) & 0 \cr 
-\sin(r) & \cos(r) & 0\cr
0 & 0 & 1
\end{pmatrix}. \ee
So here  the upper-left $2\times 2$ block is in $SO(2)$. A {\it maximal not  connected} Lie group, ${\bf K}^+\cup {\bf K}^-$, is given by the matrices 
which are either of the previous type or of the type
 \be{forma2}
K^-(r):= \begin{pmatrix} \cos(r) & -\sin(r) & 0 \cr 
-\sin(r) & -\cos(r) & 0\cr
0 & 0 & -1
\end{pmatrix}. \ee
Therefore in (\ref{forma2})   the upper-left $2\times 2$ block 
is in $O(2)$, with determinant equal to $-1$.
 We shall consider this second  case, that is, $G={\bf K}^+\cup {\bf K}^-$. 
 Remark \ref{connesso} discusses  what would change had  we chosen $G={\bf K}^+$.

\subsection{Structure of the orbit spaces $SO(3)/G$}

\begin{figure}[htb]
\centering
\includegraphics[width=0.55\textwidth, height=0.3\textheight]{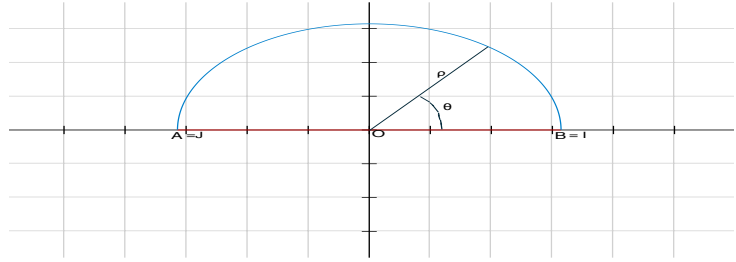}
\caption{The quotient space $SO(3)/G$.}
\label{Figso3}
\end{figure}

Following the second step of the procedure described in the previous section,  we now 
describe the structure of $M/G=SO(3)/{({\bf K}^+ \cup {\bf K}^-)}$ and its isotropy strata. We  use the Euler decomposition of $SO(3)$ from which it follows that any matrix  $X\in SO(3)$ can be written as $X=K^+(r_1)H(s)K^+(r_2)$, with $K^+(r_i)$ of the type (\ref{forma1}), and $H(s):=e^{p_1s}$, for some real $s$.  Since $K^+(r)\subset G$, $[X]=[H(s)K^+(r_2)\left(K^+(r_1)\right)^T]=[H(s)K^+(r_2-r_1)]$.
So, we can always choose as representatives of the orbits matrices of the type:
\be{classi}
H(s)K^+(r)=\begin{pmatrix} 1 & 0 & 0 \cr 
0 & \cos(s) & \sin(s) \cr
0 & -\sin(s) &  \cos(s)
\end{pmatrix} 
\begin{pmatrix} \cos(r) & \sin(r) & 0 \cr 
-\sin(r) & \cos(r) & 0\cr
0 & 0 & 1
\end{pmatrix}= \begin{pmatrix} \cos(r) & \sin(r) & 0 \cr 
-\sin(r)\cos(s) & \cos(r)\cos(s) & \sin(s) \cr
\sin(r)\sin(s) & -\cos(r)\sin(s) &  \cos(s)
\end{pmatrix} 
\ee
with $s,\, r \in [0,2\pi)$. Moreover, we have, 
\[
\begin{pmatrix} -1 & 0 & 0 \cr 
0 & -1 & 0 \cr
0 & 0 &  1
\end{pmatrix} H(s)K^+(r)\begin{pmatrix} -1 & 0 & 0 \cr 
0 & -1 & 0 \cr
0 & 0 &  1
\end{pmatrix}= \begin{pmatrix} \cos(r) & \sin(r) & 0 \cr 
-\sin(r)\cos(s) & \cos(r)\cos(s) & -\sin(s) \cr
-\sin(r)\sin(s) & \cos(r)\sin(s) &  \cos(s)
\end{pmatrix},   
\]
which changes the sign of $\sin(s)$ as compared with (\ref{classi}). 
Thus we can  assume $\sin(s) \geq 0$, so $s\in [0,\pi]$. Furthermore, we have:
\[
\begin{pmatrix} 1 & 0 & 0 \cr 
0 & -1 & 0 \cr
0 & 0 &  -1
\end{pmatrix} H(s)K^+(r)\begin{pmatrix} 1 & 0 & 0 \cr 
0 & -1 & 0 \cr
0 & 0 &  -1
\end{pmatrix}= \begin{pmatrix} \cos(r) & -\sin(r) & 0 \cr 
\sin(r)\cos(s) & \cos(r)\cos(s) & \sin(s) \cr
-\sin(r)\sin(s) &- \cos(r)\sin(s) &  \cos(s)
\end{pmatrix}, 
\]
so we can  also assume $r\in [0,\pi]$. It follows that each equivalence class has an element of the form (\ref{classi}), with $r,s\in [0,\pi]$. 
{By equating two matrices of the form 
(\ref{classi}) for different values of the pairs 
$(r,s)$,  one can see that such a correspondence is one to one  unless $s=\pi$.} In this 
case, all the matrices $H(\pi)K^+(r)$ (which give the set ${\bf K}^-$)  are equivalent. So if 
$s\in [0,\pi)$ and $r\in [0,\pi]$, each $H(s)K^+(r)$  represents 
a unique orbit, while  if $s=\pi$, since they are all equivalent,   the choice of $r$ is irrelevant. We can therefore represent  $SO(3)/G$ as the upper part of a disc of radius $\pi$, where if $\rho$ and $\theta$ are the polar coordinate, we have $\rho \in [0,\pi]$ with $\rho=\pi-s$, and $\theta\in [0, \pi]$ with $\theta=r$ (see Figure \ref{Figso3}).

\br{parametri}
If $[X_1]=[X_2]$, then $(X_1)_{3,3}=(X_2)_{3,3}$, and also 
the trace is preserved. So, from any element  $X$ 
of a given equivalence class, we can  compute the two 
parameters $s, \, r\in [0,\pi]$ of equation (\ref{classi}), by setting:
\be{calcoloparametri}
s:= \arccos (X)_{3,3},   \  \  \begin{array}{ll}
                                          r=\arccos \left( \frac{(X)_{1,1} + (X)_{2,2}}{1+X_{3,3}}\right)   &  \text{ if  $(X)_{3,3}\neq -1$},  \\
					r=0  &  \text{ if  $(X)_{3,3}= -1$} \end{array}
					\ee
From these values we have also the two values of $\rho=\pi-s$ and $\theta=r$. So there is a one to one, onto, readily  computable correspondence between points in the half 
disc in Figure \ref{Figso3} and orbits  in $SO(3)/G$.    
\er
The point $\rho=\pi$ and $\theta=0$ ($B$ in Figure 1), represents the Identity matrix, while the point $\rho=\pi$ and $\theta=\pi$ ($A$ in Figure \ref{Figso3}) gives the matrix:
\be{menoI}
J:=\begin{pmatrix} -1 & 0 & 0 \cr 
0 & -1 & 0 \cr
0 & 0 &  1
\end{pmatrix}.
\ee
Both these matrices are fixed points for the action of $G$, so they are the only matrices in their orbit, and their isotropy group is the entire group $G$. 

The points with $\rho=\pi$ and $\theta\in (0,\pi)$ give the matrices 
in ${\bf K}^+$, except for the identity ${\bf 1}$ and the matrix $J$ defined in 
(\ref{menoI}). 
The matrices in ${\bf K}^+$ commute, 
and it holds that  $K^- (v) K^+ (r)(K^- (v))^{T}=K^+ (-r)$, thus the orbits  of these elements contain two matrices, and we took as representative the one with $\sin(r)>0$.  Their isotropy group is ${\bf K}^+$.

The origin, i.e. the point with $\rho=0$ and 
$\theta=r$ arbitrary, corresponds to the matrices:
\be{origine}
\begin{pmatrix} \cos(r) & \sin(r) & 0 \cr 
\sin(r)& -\cos(r) & 0 \cr
0 &0 & -1
\end{pmatrix}.
\ee
These matrices are all equivalent, and their isotropy groups are all conjugate to:
\be{vudoppio}
W=\left\{ \begin{pmatrix} 1 & 0 & 0 \cr 
0 & -1 & 0 \cr
0 & 0 &  -1
\end{pmatrix}, \ \begin{pmatrix} -1 & 0 & 0 \cr 
0 & 1 & 0 \cr
0 & 0 &  -1
\end{pmatrix},  \ J, \ {\bf 1} \right\}
\ee
which is the isotropy group of the matrix with $r=0$.

The matrices with $\theta=\pi$ and $\rho\in (0,\pi)$, are the classes of the 
{\it symmetric matrices} in $SO(3)$. It can be seen that their isotropy group is conjugate to the one given by
\be{vu}
V=\left\{ \begin{pmatrix} 1 & 0 & 0 \cr 
0 & -1 & 0 \cr
0 & 0 &  -1
\end{pmatrix},  \ {\bf 1} \right\}
\ee

The matrices with $\theta=0$ and $\rho\in (0,\pi)$, 
correspond to   matrices in $SO(3)$ of the type:
\be{strane}
 \begin{pmatrix} a & b & c \cr 
b & d & f \cr
-c & -f&  g
\end{pmatrix},
\ee
Their isotropy group is, again, conjugate to $V$, as in the symmetric case.

The matrices which are in the interior of the half disc, have a trivial isotropy group, i.e., composed of  only  the identity matrix. This is the regular part of $SO(3)/G$ 
while the boundary of the half disc corresponds  to the singular part. 

Summarizing, the isotropy types  of $SO(3)$ are given by $(\{\bf 1\})$, $(V)$, $(W)$ in (\ref{vu}) and (\ref{vudoppio}), $({\bf K}^+)$, and $({\bf K}^+ \cup {\bf K}^-)$, with the partial ordering 
$$
(\{{\bf 1} \}) \leq (V) \leq (W) \leq ({\bf K}^+ \cup {\bf K}^-), 
$$
and 
$$(\{{\bf 1} \})  \leq ({\bf K}^+) \leq ({\bf K}^+ \cup {\bf K}^-).$$
$M_{({\bf K}^+ \cup {\bf K}^-)}$ is composed by the matrices ${\bf 1}$ and $J$, $M_{(W)}$
 are the matrices in (\ref{origine}), $M_{(V)}$ are matrices which are either symmetric or of the form (\ref{strane}), $M_{({\bf K}^+)}$ are the matrices in ${\bf K}^+$ except for ${\bf 1}$ and $J$,  $M_{(\{ {\bf 1} \})}$ are all the remaining matrices. The corresponding strata on the orbit space  
 (half disc) are indicated in Figure \ref{Figso3}.    

\subsection{Cut locus and critical locus}

We shall now apply the results given in the previous two sections to determine the cut 
locus ${\cal CL}(SO(3))$  and the critical  locus ${\cal CR} (SO(3))$. The cut locus was also described in \cite{BoscaRossi} using a different method.  Following what suggested in Remark \ref{sugge3}, we analyze the singular points, first.

\bp{prima}
All the matrices that correspond to $\rho=\pi$ and $\theta\in (0,\pi]$ (these are 
all the matrices in ${\bf K}^+$ except the ${\bf 1}$)  are in ${\cal CL}(SO(3))$, and so also in the 
${\cal CR}(SO(3))$ (cf. Proposition \ref{Propo1}). 
\ep
\bpr
Fix a matrix $X_f\in {\bf K}^+$ and let $A_p=\alpha p_1+\beta p_2$, be the matrix giving  the minimizing geodesic that appear in  equation (\ref{geodesicKP}) for $X_f$.
If this matrix is not in ${\cal CL}(SO(3))$, then, using Proposition \ref{ProposizioneFrancesca}, it must hold:
\[
[A_p ,K^+ ]=0,
\]
for all $K^{+}\in {\bf K}^+$, since ${\bf K}^+$ is contained in the isotropy group (indeed ${\bf K}^+$ {\it is} the isotropy group for all values of $\theta\in (0,\pi)$, while for $\theta=\pi$ the isotropy group is all $G$). The   previous equality holds for all  $K^+$ if and only if $A_p=0$, which is not possible since $||A_p||=1$.
So $X$ is in  the cut locus, and also in the critical locus.
\epr

The next proposition proves that all the symmetric matrices (which correspond to the segment $O-A$ in Figure \ref{Figso3}) are in the cut locus.
\bp{seconda}
The matrices corresponding to $\rho=0$ and to $\rho\in (0,\pi)$ and $\theta=\pi$ (these are the matrices which correspond to the origin and to the segment $(A,O)$ in the Figure \ref{Figso3}) are in ${\cal CL}(SO(3))$, and so also in ${\cal CR}(SO(3))$. 
\ep
\bpr
Fix a symmetric matrix $X_f$. Its isotropy group is conjugate either to $W$ in (\ref{vudoppio}) (if $\rho=0$) or to $V$ in (\ref{vu}).  By continuity  the geodesic  from {\bf 1} to $X_f$ must contain matrices whose isotropy group is different from the one of $X_f$, so by using Corollary 
\ref{coro77} we get that $X_f$ lies in the cut locus.
\epr

Now we will prove that all the remaining matrices, i.e. the ones corresponding to the 
open segment $(OB)$ and the regular part (the interior of the disc) are neither on the ${\cal CL}(SO(3))$ nor in the critical locus ${\cal CR}(SO(3))$.

We know, that the geodesic are analytic curves given by equation (\ref{geodesicKP}). Here we may choose as $\mathcal{A}=\text{span }\{p_1\}$, thus the geodesic are given by:\footnote{Here we use the calculation of \cite{BoscaRossi} section 3.2.1.}
\be{traiettoria}
[X(t)]=[e^{-\alpha k t}e^{(\alpha k +p_1)t}]=
\ee
\[
\left[ \begin{pmatrix} \frac{1+C_1\alpha^2}{1+\alpha^2} \cos(\alpha t)+ C_2\alpha \sin(\alpha t) &C_1\sin(\alpha t) -C_2\alpha \cos(\alpha t) & C_3\cos(\alpha t)-C_2 \sin(\alpha t)  \cr 
-\frac{1+C_1\alpha^2}{1+\alpha^2} \sin(\alpha t)+ C_2\alpha \cos(\alpha t)  & C_1\cos(\alpha t) +C_2\alpha \sin(\alpha t) & -C_2\cos(\alpha t)-C_3\sin(\alpha t)  \cr
C_3& C_2 &  \frac{C_1+\alpha^2}{1+\alpha^2}
\end{pmatrix}\right],
\]
where 
\[
C_1=\cos\left( \sqrt{(1+\alpha^2)}t\right), \  \   C_2= \frac{\sin\left( \sqrt{(1+\alpha^2)}t\right)}{ \sqrt{1+\alpha^2}}, \  \  C_3=\frac{\alpha\left( 1- \cos\left( \sqrt{(1+\alpha^2)}t\right)\right)}{1+\alpha^2}.
\]

The next proposition gives the optimal time to reach the matrices with $\rho=0$, i.e., the ones corresponding to the origin of the half disc as in (\ref{origine}). 
\bp{origine2}
The optimal geodesic to reach any $X_f$ such that $[X_f]=\left[\begin{pmatrix} 1 & 0 & 0 \cr 
0 & -1 & 0 \cr
0 & 0 &  -1
\end{pmatrix}\right]$  must have the parameter $\alpha$ of equation (\ref{traiettoria}) equal to $0$, and the minimum time to reach $X_f$ is $\pi$.
\ep
\bpr
Since the conjugation by elements of $G$ does not change the $3,3$ element, letting $T$ the minimum  time to reach $X_f$, we must have (see equation  (\ref{traiettoria})):
\[
\frac{\cos\left( \sqrt{(1+\alpha^2)}T\right)+\alpha^2}{1+\alpha^2}=-1.
\]
The previous equality can hold if and only if $\alpha=0$. Moreover we must have $\cos(T)=-1$. 
Thus the minimum time $T$ is equal to $\pi$.
\epr

The next proposition proves that the matrices  in the singular part which correspond to  the segment $O-B$ in Figure \ref{Figso3}, are neither on the cut locus nor on the critical locus.
In particular this implies that the projection of the geodesics reaching these matrices lies  all in the segment, since each point of these trajectories has to have the same isotropy group. 

\bp{terza}
Fix the matrix $X_f$ that corresponds to $\theta=0$ and $\rho=\pi-s$, with $s\in(0,\pi)$ as in (\ref{strane}). Then this matrix is not on the cut locus nor on the critical locus, and the minimum  time $T$ to reach $X_f$ from ${\bf 1}$ is $T=s$.
\ep
\bpr
Fix a matrix $X_f$ that corresponds to $\theta=0$ and $\rho=\pi-s$, with $s\in(0,\pi)$, i.e. such that 
$$[X_f]=\left[\begin{pmatrix} 1 & 0 & 0 \cr 
0 &\cos(s) & \sin(s) \cr
0 & -\sin(s) & \cos(s)
\end{pmatrix}\right].$$
These are matrices of the form (\ref{strane}). 
First we prove that necessarily the geodesic reaching $X_f$ must have 
$\alpha=0$. Let $\gamma(t)$ be a geodesic  with $\alpha=0$.  Then, by Proposition \ref{origine2}, its projection is optimal until $t=\pi$, thus $\gamma(t)$ is optimal until $t=\pi$.
Moreover since its projection at time  $t=s$ is equal to $H(s):=e^{p_1s}$, we have $\gamma(s)=X_f$, and $s$ is the minimum time, since the minimum time is the same for equivalent matrices (cf. Proposition \ref{Propo1b}). 
If there was another trajectory reaching optimally $X_f$, with $\alpha\neq 0$, and call this trajectory $\tilde{\gamma}(t)$, then the trajectory:
\[
\eta(t)=\left\{ \begin{array}{ll}
\tilde{\gamma}(t)   & t\in [0,s)\\
\gamma(t)          &  t\in [s,\pi],
\end{array}\right.
\]
would also be an optimal trajectory to the origin, which contradicts the fact that all geodesics are analytic.

Assume now that $X_f$ is on the cut locus. Then there exist two optimal trajectories both with $\alpha=0$, so  $\gamma_i(t)=e^{A_{p}^it}$ such that, 
\[
X_f=e^{A_p^1s}=e^{A_p^2s}.
\]
Since every two  Abelian subagebras in $\mathcal{P}$ are conjugate by an element of ${\bf K}^+$, there must exist a matrix $K^+ \in {\bf K}^+$ such that
\[
\text{span }\{ A_p^2\}=K^+ \text{span }\{ A_p^1\} (K^+)^{T},
\]
However,  since these spans are one dimensional, we must have 
\be{piuomeno}
A_p^2=\pm K^+A_p^1(K^+)^{T}. 
\ee
Thus
\[
X_f=e^{A_p^2s}=\left\{
\begin{array}{ll}
K^+e^{A_p^1s}(K^+)^{T}  &  \text{ if   (\ref{piuomeno}) is verified with   $+1$} \\
K^+ e^{-A_p^1s}(K^+)^T&  \text{ if   (\ref{piuomeno}) is verified with   $-1$} 
\end{array} \right.
\]
In the first case, we have that $K^+$ must be in the isotropy group of $X_f$. On the other hand the isotropy group of $X_f$ is conjugate to the group $V$ of equation (\ref{vu}), thus it contains two elements, one is the identity and the other must have $-1$ in the $3,3$ position. Thus necessarily since $K^+$ has $+1$ in the $3,3$ position, we must have $K^+={\bf 1}$, and so $A_p^2=A_p^1$.
 In the second case, $X_f$ is conjugate via an element of ${\bf K}^+$ to $X_f^{-1}=X_f^T$. Writing the third column of the relation $X_f K^+=K^+ X_f^T$ using the formula (\ref{strane}) with $K^+:=\begin{pmatrix} K_1^+ & 0 \cr 0 & 1\end{pmatrix}$ as 
 \be{pppl}
 \begin{pmatrix}
 c \cr f \cr g
 \end{pmatrix}=\begin{pmatrix}
-K_1^+ \begin{pmatrix}  c \cr f\end{pmatrix} \cr g  
 \end{pmatrix},
 \ee
we have that the $2 \times 2$ matrix $K_1^+ \in SO(2)$ has an eigenvalue in $-1$ (unless $c$ and $f$ are equal to zero which is to be excluded since $g\not= \pm 1$). Therefore 
\[
K^+=\begin{pmatrix}-1 & 0 & 0 \cr 
0 & -1 & 0 \cr 
0 & 0 & 1 \end{pmatrix}
.\] 
Using this in $A_p^2=-K^+ A_p^1 (K^+)^T$ and the general expression for $A_p^1$, we find again $A_p^2=A_p^1$.  

Therefore $X_f$ is not on the cut locus. Moreover, since the projection of the trajectory is optimal until $t=\pi>s$, the matrix $X_f$ is not on the critical locus either.
\epr

\subsection{The optimal synthesis}

The last proposition has characterized the minimizing geodesics for points corresponding to the interval $O-B$ in Figure \ref{Figso3}, while Proposition \ref{origine2} has 
given the minimizing geodesic and optimal time for points corresponding to the origin, i.e. matrices in $\bf{K}^-$,  in Figure \ref{Figso3}.
We now consider the geodesics leading to the remaining pieces of the singular part of $SO(3)/G$. Then we put all things together to describe the full optimal synthesis. 

The geodesic curves  given in equation (\ref{traiettoria}) depend on  the parameter $\alpha$ which varies in $\RR$. However 
both parameters $\rho$ and $\theta$ which characterize the points of the equivalence classes in the orbit space are even function of $\alpha$ (see equation (\ref{calcoloparametri})), so in the analysis in the orbit space, we can restrict ourselves to values  $\alpha\geq 0$.

The next Proposition provides the optimal time to reach any matrix with $\rho=\pi$, i.e.,  all the matrices in ${\bf K}^+$.

\bp{bordo}
Assume $X_f\in {\bf K}^+$, then $[X_f]=\{X_f,X_f^T\}$, and let $\theta \in (0,\pi]$ be the value of the parameter of equation (\ref{classi}), which together with $\rho=\pi$ gives 
 the equivalence class $[X_f]$. Then the minimum time $T$ to reach $X_f$ is given by
\[
T=\sqrt{\theta(4\pi-\theta)},
\]
and the optimal value of the parameter $\alpha$  to reach  $[X_f]$ is $\alpha=\frac{2\pi-\theta}{\sqrt{\theta(4\pi-\theta)}}$.
\ep
\bpr
First notice that necessarily $\alpha\neq 0$, since all the trajectories corresponding to $\alpha=0$ have $\theta=0$.
Since the equivalence class of $X_f$ consists of only two elements (which coincide when $\theta=\pi$) and these elements have $0$ in the $3,1$ and $3,2$ position, and $1$ in the $3,3$ position, for $t=T$ we must have in equation (\ref{traiettoria}), $C_2=C_3=0$ and $\frac{C_1+\alpha^2}{1+\alpha^2}=1$, which implies:
\[
C_2=\sin\left(\sqrt{(1+\alpha^2)} T\right)=0 \ \text{ and } \ C_1=\cos\left(\sqrt{(1+\alpha^2)} T\right)=1,
\]
thus we must have:
\be{bordouno}
\sqrt{(1+\alpha^2)} T= 2m\pi,
\ee
for some $m\in \NN$. Moreover at time $T$, we have:
\[
[X(T)]= \left[ \begin{pmatrix} \cos(\alpha T) & \sin(\alpha T) & 0 \cr 
-\sin(\alpha T) & \cos(\alpha T) & 0\cr
0 & 0 & 1
\end{pmatrix}\right],
\]
which implies 
\be{bordodue}
\cos(\alpha T)=\cos(\theta) \  \  \  \Rightarrow  \  \   \  \alpha T=\pm \theta + 2p \pi,
\ee
for some $p\in \NN$. We will treat the $\pm \theta$ sign separately.

{\em{ Case +1}} Assume that equation (\ref{bordodue}) holds with the $+1$ sign. Since $\sqrt{(1+\alpha^2)} T>\alpha T$, we must have $p\leq m-1$. From equation (\ref{bordouno}) and (\ref{bordodue}), we have:
\[
\frac{2m\pi}{\sqrt{(1+\alpha^2)}}=\frac{\theta + 2p \pi}{\alpha}.
\]
The previous equality implies:
\[
\alpha= \frac{\theta + 2p \pi}{\sqrt{(4m^2\pi^2-({\theta + 2p \pi})^2)}},
\]
and consequently:
\[
T_{m,p}={\sqrt{(4m^2\pi^2-({\theta + 2p \pi})^2)}}.
\]
The value of $T_{m,p}$, for each fixed $m$,  is minimum when $p$ is maximum, i.e. $p=m-1$, and its  minimum value is 
\[
T_{m,m-1}=\sqrt{(2\pi-\theta)(4m\pi+\theta-2\pi)}, 
\]
which is minimum when $m=1$ and we have $T_{1,0}:=T_{1,0}^+=\sqrt{(2\pi-\theta)(2\pi+\theta)}$.

{\em{ Case -1}} Assume that equation (\ref{bordodue}) holds with the $-1$ sign. Imposing again $\sqrt{(1+\alpha^2)} T>\alpha T$, we now get  $p\leq m$. From equation (\ref{bordouno}) and (\ref{bordodue}) we have:
\[
\frac{2m\pi}{\sqrt{(1+\alpha^2)}}=\frac{-\theta + 2p \pi}{\alpha}.
\]
The previous equality implies:
\[
\alpha= \frac{-\theta + 2p \pi}{\sqrt{(4m^2\pi^2-({-\theta + 2p \pi})^2)}},
\]
and consequently:
\[
T:=T_{m,p}={\sqrt{(4m^2\pi^2-({ 2p \pi-\theta })^2)}}.
\]
Again $T^{-}_{m,p}$, for each fixed $m$ is minimum when $p$ is maximum. Therefore,  we now  take $p=m$, and we get:
\[
T^{-}_{m,m}=\sqrt{\theta(4m\pi-\theta)}, 
\]
which is minimum when $m=1$ and we have $T^{-}_{1,1}=\sqrt{\theta(4\pi-\theta)}$.

Since $\theta\leq \pi$, we have $T^{-}_{1,1}\leq T^{+}_{1,0}$, thus the minimum time is $T=\sqrt{\theta(4\pi-\theta)}$ with  the corresponding  $\alpha=\frac{2\pi-\theta}{\sqrt{\theta(4\pi-\theta)}}$.
\epr

From the previous Proposition, since $\theta\in (0,\pi)$, we have that for $\alpha\geq \frac{1}{\sqrt 3}$, all the geodesics are optimal until time 
$T=\frac{2\pi}{\sqrt{1+\alpha^2}}$, when they reach the boundary of the disc. It is clear that $T$ is an increasing function of $\alpha$, with maximum equal to $\pi\sqrt 3$, which corresponds to the trajectory reaching  the matrix $J$.  The trajectory corresponding to $\alpha=0$ lies on the segment  $(O,B)$ and it is optimal until time $T=\pi$, when it reaches the origin. The trajectories corresponding to $\alpha\in (0,\frac{1}{\sqrt 3})$ are 
optimal until they reach the segment $(A,O)$, which correspond to the symmetric matrices.
We know from  Proposition \ref{seconda}  that these matrices are on the cut locus. For a given $\alpha$, the time $T$ where the corresponding geodesic loses optimality, can be numerally estimated, and it is always between $\pi$ and $\sqrt 3\pi$.

Thus all elements are reached in time $T\leq \sqrt 3\pi$.
See Figure 2 for the shape of the optimal trajectories, the red curve is the optimal curve with $\alpha=\frac{1}{\sqrt 3}$, the black curves correspond to bigger values of $\alpha$ and loose optimality at the boundary of the circle, while the blue curves correspond to smaller values of $\alpha$ and loose optimality at the segment $(A,O)$.

\begin{figure}[htb]
\centering
\includegraphics[width=0.7\textwidth, height=0.4\textheight]{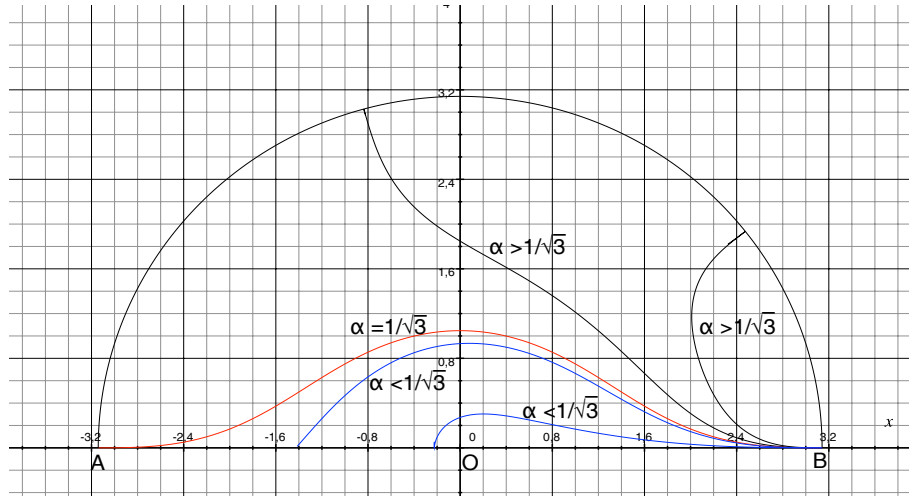}
\caption{Optimal trajectories}
\label{geodesichefigura1}
\end{figure}
Figure 3 describes the optimal synthesis according to the third step of the procedure given in the previous section, that is, it gives the boundaries of the 
 reachable sets at any time $t$. To draw these curves, for a given time $T$ one finds the values of $\alpha$ such that the corresponding trajectory at time $T$ lies on the boundary, and these are parametric curves with $\alpha$ as a parameter in the given interval. For $T<\pi$, the boundary is given varying $\alpha$ from $0$, until the boundary of the circle is reached, for $T>\pi$, the parameter $\alpha$ has to be chosen from the values that correspond to the segment $(A,O)$ until it again reaches the boundary of the circle. So the behavior changes at the curve in red corresponding to $T=\pi$. 
\begin{figure}[htb]
\centering
\includegraphics[width=0.6\textwidth, height=0.3\textheight]{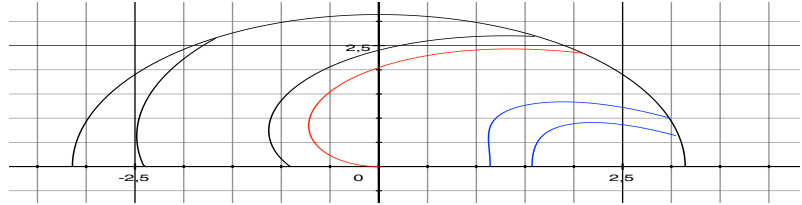}
\caption{Reachable Sets}
\label{geodesichefigura2}
\end{figure}

\br{connesso}
To derive all the previous results we have taken as symmetry group $G=K^+\cup K^-$. We could have done a similar  analysis, taken as a group of symmetries only the connected component containing the origin, i.e. $\tilde{G}=K^+$. In this case as representatives of equivalent classes we could take again matrices of the type 
(\ref{classi}), but now, while $s\in [0,\pi]$, we may allow $r\in(-\pi,\pi]$. So the quotient space turns out to be the all disk of radius $\pi$, instead of only the upper part.
Here the boundary, represents the matrices in $K^+$, that now are all fix points and the center are the matrices in $K^-$, which are again all equivalent.
It is easy to see that this two sets give the singular part of $SO(3)/(K^+)$, while the interior of the disk is all in the regular part.
The trajectories in the quotient space are given by the trajectories we have found previously and the one that are the symmetric with respect to the $x-$axis, 
this can be easily seen, since the two parameters $s,\, r$ can be found using, as before, equations (\ref{parametri}), but while $s$ is the same, for $r$ we have two choices, the $r$  given in (\ref{parametri}) and its opposite (see also figure 4).
\er

\begin{figure}[htb]
\centering
\includegraphics[width=0.8\textwidth, height=0.4\textheight]{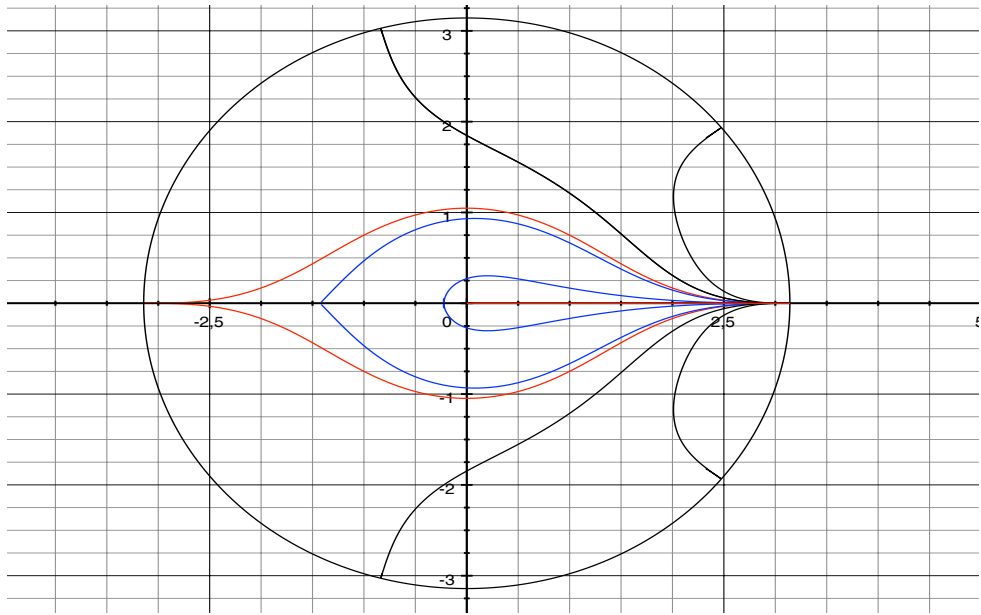}
\caption{Geodesics when $G=K^+$}
\label{geodesichefigura3}
\end{figure}



\section*{Acknowledgement} Domenico D'Alessandro's  research was  supported by ARO MURI grant W911NF-11-1-0268. Domenico D'Alessandro also would like to thank the Institute of Mathematics and its Applications in Minneapolis and the Department of Mathematics at the University of Padova, Italy, for kind hospitality during part of this work. 


\begin{thebibliography}{99}


\bibitem{ABB} A. Agrachev, D. Barilari and U. Boscain, Introduction to Riemannian and sub-Riemannian geometry, Lecture Notes SISSA, Trieste, Italy, 2011. 

\bibitem{AS} A. Agrachev and Y. Sachkov, {\it Control Theory from the Geometric Viewpoint}, Encyclopaedia of Mathematical Sciences, {\bf 87}, 2004, Springer-Verlag Berlin-Heidelberg.


\bibitem{Noi} F. Albertini and D. D'Alessandro, Minimum time optimal synthesis for two  level quantum systems, {\it Journal of Mathematical Physics}, {\bf 56}, 012106 (2015).


\bibitem{OptimalityModels} F. Albertini and D. D'Alessandro, Time Optimal Simultaneous Control of Two Level Quantum Systems, submitted to {\it Automatica}.

\bibitem{Dimitry} D. Alekseevsky, A. Kriegl, M. Losik and P. W. Michor, The Riemannian geometry of orbit spaces. The metric, geodesics, and integrable systems, {\it Publ. Math. Debrecen}, {\bf 62} (2003), 247-276. 




\bibitem{OptRes} U. Boscain, T. Chambrion, and J.P. Gauthier, On the K+P problem for a three-level quantum system: Optimality implies resonance, {\it Journal of Dynamical and Control Systems}, Vol. 8, No. 4, October 2002, 547-572.




\bibitem{BoscaRossi} U. Boscain and F. Rossi, Invariant Carnot-Caratheodory metric on 
$S^3$, $SO(3)$ and $SL(2)$ and Lens Spaces, {\it SIAM Journal on Control and Optimization}, Vol. 47, pp. 1851-1878, (2008).

\bibitem{LTG} G. E. Bredon, {\it Introduction to Compact Transformation Groups}, Pure and Applied Mathematics, Vol. 46,  Academic Press, New York, 1972.






\bibitem{NoiSIAM} D. D'Alessandro, F. Albertini  and R. Romano, Exact algebraic conditions for indirect controllability of quantum systems, {\it SIAM Journal on Control and Optimization}, 2015 53:3, 1509-1542. 








\bibitem{dodici} A. Echeverr\`ia-Enriquez, J. Mar\`in-Solano, M.C. Mun\~oz Lecanda and N. Roman-Roy, Geometric reduction in optimal control theory with symmetries, {\it Rep. Math. Phys.}, 52 (2003), pp. 89-113.




\bibitem{Bredon} G. E. Bredon, {\it Introduction to Compact Transformation Groups}, Pure and Applied Mathematics, Vol. 46,  Academic Press, New York, 1972.

\bibitem{Filippov} A. F. Filippov, On certain questions in the theory of optimal control, {\it SIAM J. on Control,} Vol 1, pp/ 78-84, 1962. 

\bibitem{quattordici} J. Grizzle and S. Markus, The structure of nonlinear control systems possessing symmetries, {\it IEEE Trans. Automat. Control}, 30, (1985), pp. 248-258.

\bibitem{quindici} J. Grizzle and S. Markus, Optimal control of systems possessing symmetries, {\it IEEE Trans. Automat. Control}, 29 (1984), pp. 1037-1040.


\bibitem{sedici} A. Ibort, T. R. De la Pe\"na, and R. Salmoni, Dirac structures and reduction of optimal control problems with symmetries, preprint 2010.

\bibitem{Jacquet} S. Jacquet, Regularity of the sub-Riemannian distance and cut locus, in {\it Nonlinear Control in the Year 2000}, Lecture Notes in Control and Information Sciences, Vol. 258 (2007), pp. 521-533. 



\bibitem{Knapp} A. Knapp, {\it Lie Groups Beyond and Introduction}, Progress in Mathematics, Vol. 140, Birkh\"{a}user Boston, 1996. 

\bibitem{venti} W. S. Koon and J. E. Marsden, The Hamiltonian and Lagrangian approaches to the dynamics of nonholonomic systems, {\it Rep. Math. Phys.}, 40 (1997), pp. 21-62.



\bibitem{Marsden1} J.E. Marsden and T.S. Ratiu, {\it Introduction to Mechanics and Symmetry}, Springer, New York, 1999.

\bibitem{Marsden2} J. E. Marsden and A. Weinstein, Reduction of symplectic manifolds with symmetry, {\it Rep. Math. Phys.} 5 (1974), pp. 121-130.

\bibitem{Martinez} E. Martinez, Reduction in optimal control theory, {\it Rep. Math. Phys.}, vol 53 (2004), No. 1, pp. 79-90.

\bibitem{EM} E. Meinrenken, Group Actions on Manifolds (lecture notes), University of Toronto, 2003.  

\bibitem{Michor1} P. Michor, {\it Isometric Actions of Lie Groups and Invariants}, Lecture Course at the University of Vienna, 1996-1997. 

\bibitem{Montgomery} R. Montgomery, {\it A Tour of sub-Riemannian geometries, their Geodesics and Applications}, volume 91 of Mathematical Surveys and Monographs, American Mathematical Society, RI, 2002.  

\bibitem{MY} R. Monti, The regularity problem for sub-Riemannian geodesics, in {\it Geometric Control and Sub-Riemannian Geometry}, G. Stefani, U. Boscain, J-P. Gauthier, A. Sarychev, and M. Sigalotti Eds, Springer INdAM Series, Volume 5 2014, pp. 313-332. 


\bibitem{trentasei} H. Nijmeijer and A. Van der Schaft, { Controlled invariance for nonlinear systems}, {\it IEEE Trans. Automat. Control,} 27, (1982), pp. 904-914.



\bibitem{Tomizu} T. Ohsawa, Symmetry reduction of optimal control systems and principal connections, {\it SIAM J. Control Optim.}, Vol. 51, No. 1, pp 96-120, (2013).








\end{thebibliography}
\end{document}